\documentclass[11pt]{article}
\usepackage{srcltx}
\usepackage{epsfig}
\usepackage{url}
\def\TV{{\hbox{\rm TV}}}
\def\bZ{{\mathbf{Z}}}
\def\CTV{{{\cal T}{\cal V}}}
\def\cE{{\cal E}}
\def\tv{{\hbox{\scriptsize TV}}}
\def\DFT{{\hbox{DFT}}}
\def\cC{{\cal C}}
\def\dft{{\hbox{\scriptsize\rm DFT}}}
\def\CG{CndG}
\usepackage{psfrag}
\usepackage{graphicx}
\usepackage[latin1]{inputenc}
\usepackage{subfigure}
\usepackage{fullpage}
\usepackage{hyperref}
\usepackage{mathtools}
\usepackage{amsfonts,amsmath,amssymb,amsthm}
\usepackage{color}
\usepackage{booktabs}

\usepackage{float,enumitem}
\setlist{itemsep=3pt,parsep=0pt,topsep=3pt}

\floatstyle{ruled}
\newfloat{algorithm}{tbp}{loa}
\floatname{algorithm}{Algorithm}

\newtheorem{thm}{Theorem}

\newtheorem{lemma}[thm]{Lemma}

\def\qed{\ \hfill$\square$\par\smallskip}

\def\Conv{\hbox{\rm Conv}}

\def\Row{\hbox{\rm Row}}

\def\Tr{{\mathop{\hbox{\rm Tr}}}}
\def\Opt{{\mathop{\hbox{Opt}}}}

\newcommand{\cN}{I\!\! N}

\newcommand{\C}{{\cal C}}

\newcommand{\half}{ \mbox{\small$\frac{1}{2}$}}

\def\Ker {{\rm Ker}}

\newcommand{\be}{\begin{eqnarray}}
\newcommand{\ee}[1]{\label{#1}\end{eqnarray}}

\newcommand{\ese}{\end{eqnarray*}}
\newcommand{\bse}{\begin{eqnarray*}}
\newcommand{\rf}[1]{~(\ref{#1})}

\newtheorem{proposition}{Proposition}
\newtheorem{theorem}{Theorem}

\def\Argmin{\mathop{\hbox{\rm Argmin$\,$}}}

\def\Argmax{\mathop{\hbox{\rm Argmax$\,$}}}

\def\cA{{\cal A}}

\def\cN{{\cal N}}
\def\bS{{\mathbf{S}}}
\def\bR{{\mathbf{R}}}

\title{Conditional Gradient Algorithms\\ for Norm-Regularized Smooth Convex Optimization}
\date{}
\author{Zaid Harchaoui
\thanks{LJK,
INRIA Rh\^one-Alpes,  655 Avenue de l'Europe, Montbonnot, 38334 Saint-Ismier
France	
 {\tt zaid.harchaoui@inria.fr}}
\and
Anatoli Juditsky\thanks{LJK,
Universit\'e J. Fourier, B.P. 53, 38041 Grenoble
Cedex 9, France, {\tt anatoli.juditsky@imag.fr}}
\and Arkadi Nemirovski\thanks{Georgia Institute
 of Technology, Atlanta, Georgia
30332, USA, {\tt nemirovs@isye.gatech.edu}\newline
Research of
the third author was supported by the ONR grant N000140811104 and NSF grants DMS 0914785, CMMI 1232623}
}

\begin{document}
\maketitle
\vspace{-0.7cm}

\begin{abstract}
Motivated by some applications in signal processing and machine learning, we
consider two convex optimization problems where, given a cone $K$, a norm
$\|\cdot\|$ and a smooth convex function $f$, we want either 1) to minimize the
norm over the intersection of the cone and a level set of $f$, or 2) to minimize
over the cone the sum of $f$ and a multiple of the norm. We focus on the case
where (a) the dimension of the problem is too large to allow for interior point
algorithms, (b) $\|\cdot\|$ is ``too complicated'' to allow for computationally
cheap Bregman projections required in the first-order proximal gradient algorithms. On
the other hand, we assume that {it is relatively easy to minimize linear forms over the intersection of $K$ and the unit $\|\cdot\|$-ball}. Motivating examples are given by the nuclear norm with $K$
being the entire space of matrices, or the positive semidefinite cone in the
space of symmetric matrices, and the Total Variation norm on the space of 2D
images. We discuss versions of the Conditional Gradient algorithm
capable to handle our problems of interest,
provide the related theoretical efficiency estimates and outline some
applications.
\end{abstract}

%
%
\section{Introduction}
We consider two {\sl norm-regularized} convex optimization problems as  follows:
\be
\label{eq:starteq0}
\mbox{[norm minimization]}&&\min_{x\in K}\,\|x\|,\;\;\mbox{subject
to}\;\;f(x)\le \delta,\\
\mbox{[penalized minimization]}&&   \min_{x  \in K } \,
    f(x) + \kappa \Vert x \Vert
\ee{eq:starteq1}
where $f$ is a convex function with Lipschitz continuous gradient, $K$ is a
closed convex cone in a Euclidean space $E$, $\|\cdot\|$ is some norm,
$\delta$ and $\kappa$ are positive parameters. Problems such
as  such as (\ref{eq:starteq0})
and (\ref{eq:starteq1}) are of
 definite interest for signal processing and machine learning. In these
applications, $f(x)$ quantifies the  discrepancy between the observed noisy
output of some parametric model and the output of the model with candidate
vector $x$ of parameters. Most notably, $f$ is the quadratic penalty:
$f(x)=\frac{1}{2}\|\cA x - y\|_2^2$, where $\cA x$ is the ``true'' output of the
linear regression model $x\mapsto \cA x$, and $y=\cA x_*+\xi$, where $x_*$ is
the vector of true parameters, $\xi$ is the observation error, and $\delta$ is
an a priori upper bound on ${1\over2}\|\xi\|_2^2$. The cone $K$ sums up a priori
information on the parameter vectors (e.g., $K=E$ -- no a priori information at
all, or $E=\bR^p,$ $K=\bR^p_+$, or $E=\bS^p$, the space of symmetric $p\times
p$ matrices, and $K=\bS^p_+$, the cone  of positive semidefinite matrices, as
is the case of covariance matrices recovery). Finally, $\|\cdot\|$ is a
regularizing norm ``promoting''  a desired property of the recovery, e.g., the
sparsity-promoting norm $\ell_1$ on $E=\bR^n$, or the low rank promoting nuclear
norm on $E=\bR^{p\times q}$, or the Total Variation (TV) norm, as in image
reconstruction.\par
 In the large-scale case, first-order algorithms of proximal-gradient type are popular
to tackle such problems,
see~\cite{Sra:Nowozin:Wright:2010} for a recent overview.
Among them, the celebrated Nesterov optimal gradient methods for smooth and
composite minimization~\cite{Nesterov:2004,Nesterov:2007,Nes:Nem:2013},
and their stochastic approximation counterparts~\cite{Lan:2012}, are now
state-of-the-art in compressive sensing and machine learning.
These algorithms enjoy the best known so far theoretical estimates
 (and in some cases, these estimates are the best
possible for the first-order algorithms).
For instance, Nesterov's algorithm for penalized minimization
\cite{Nesterov:2007,Nes:Nem:2013} solves~\rf{eq:starteq1} to accuracy $\epsilon$
in $O(D_0\sqrt{L/\epsilon})$ iterations,
where $L$ is the properly defined Lipschitz constant of the gradient of $f$,
and $D_0$ is the initial distance to the optimal set, measured in the norm
$\Vert \cdot
\Vert$.
However, applicability and efficiency of proximal-gradient algorithms in the
large-scale case require from the problem to possess
``favorable geometry'' (for details, see \cite[Section A.6]{Nes:Nem:2013}). To
be more specific, consider proximal-gradient algorithm for convex minimization problems
of the form
\begin{equation}
\label{eq:starteq3}
\min_{x} \; \{f(x):\;\|x\|\le 1, \,x\in K\} \, .
\end{equation}
The comments to follow, with
slight modifications, are applicable to problems such as (\ref{eq:starteq0})
and (\ref{eq:starteq1}) as well.
In
this case, a proximal-gradient algorithm operates with a ``distance generating
function'' (d.g.f.) defined on the domain of the problem and $1$-strongly convex
w.r.t. the norm $\Vert\cdot\Vert$. Each step of the  algorithm requires
minimizing the sum of the d.g.f. and a linear form. The efficiency estimate of
the algorithm depends on the variation of the d.g.f. on the domain and on
regularity of $f$ w.r.t. $\|\cdot\|$ \footnote{i.e., the Lipschitz constant of
$f$ w.r.t. $\Vert\cdot\Vert$ in the nonsmooth case, or the Lipschitz constant of
the gradient mapping $x\mapsto f'(x)$ w.r.t. the norm $\Vert\cdot\Vert$ on the
argument and the conjugate of this norm on the image spaces in the smooth
case.}. As a result, in order for a proximal-gradient algorithm to be practical in the
large scale case, two ``favorable geometry'' conditions should be met: (a)  the
outlined sub-problems should be easy to solve, and (b)  the variation of the
d.g.f. on the domain of the problem should grow slowly (if at all) with
problem's dimension. Both these conditions indeed are met in many applications;
see, e.g., \cite{Bach:family:2012,Jud:Nem:2012} for examples.
This explains the recent  popularity of this family of algorithms.

However, sometimes conditions (a) and/or (b) are violated, and application of proximal
algorithms becomes questionable. For example, for the case of $K=E$, (b) is violated for
the usual $\|\cdot\|_\infty$-norm on $\bR^p$ or, more generally, for
$\|\cdot\|_{2,1}$ norm on the space of $p\times q$ matrices given by
\[
\|x\|_{2,1}=\max_{1\le j\le p} \|\Row_j(x)\|_2,
 \]
where $\Row_j^T(x)$ denotes the $j$-th row of $x$.
Here the variation of (any) d.g.f. on problem's domain is at least $p$.
As a result, in the case in question the theoretical iteration complexity of a
proximal algorithm grows rapidly with the dimension $p$.
Furthermore, for some high-dimensional problems which do satisfy (b), solving the
sub-problem can be computationally
challenging. Examples of such problems include nuclear-norm-based
 matrix completion, Total Variation-based image reconstruction, and
multi-task learning with a large
number of tasks and features. This corresponds to
$\Vert\cdot\Vert$ in \rf{eq:starteq0} or \rf{eq:starteq1} being the nuclear
norm~\cite{hadopaduma:2012} or the TV-norm.
\par
These limitations recently motivated alternative approaches, which do not rely
upon favorable geometry of the problem domain and/or do not
require to solve hard sub-problems at each iteration, and triggered
a renewed interest
in the Conditional Gradient ({\CG}) algorithm. This algorithm,
{also known as the Frank-Wolfe algorithm,}
which is historically
the first method for smooth constrained convex optimization, originates from
\cite{Frank1956Algorithm}, and was extensively studied in the 70-s (see,
e.g.,~\cite{Dem:Rub:1970,Dunn78,Pshe:1994} and references therein). {\CG}
algorithms
work by minimizing a linear form on the problem domain at each iteration; this
auxiliary problem clearly is easier, and in many cases -- significantly easier
than the auxiliary problem arising  in proximal-gradient algorithms.
Conditional gradient
algorithms for collaborative filtering were studied
recently~\cite{Jaggi:2010,Jaggi:2013}, some  variants and extensions
were studied in~\cite{dhm:2012,ShalevShwartzGoSh11,hadopaduma:2012}.
Those works consider constrained formulations of machine learning or signal processing problems, i.e.,
minimizing the discrepancy $f(x)$ under a constraint on the norm of the solution, as in~\rf{eq:starteq3}.
On the other hand, {\CG} algorithms for
other learning formulations, such as norm minimization \rf{eq:starteq0}
or penalized minimization \rf{eq:starteq1} remain open issues.
An exception is the work of~\cite{dhm:2012,hadopaduma:2012},
where a Conditional Gradient algorithm for penalized minimization was studied,
although the efficiency estimates obtained in that paper were suboptimal.
In this paper, we present {\CG}-type algorithms aimed at solving norm minimization
and penalized norm minimization problems and provide theoretical efficiency
guarantees for these algorithms.\par
The main body of the paper is organized as follows. In Section \ref{sec:prb},
we present detailed setting of problems (\ref{eq:starteq0}), (\ref{eq:starteq1})
along with basic assumptions on the ``computational environment'' required by
the {\CG}-based algorithms we are developing. These algorithms and their efficiency
bounds are presented in Sections \ref{sec:po} (problem (\ref{eq:starteq0})) and
\ref{sec:cocg} (problem (\ref{eq:starteq1}). In Section \ref{sec:examples} we
outline some applications, and in Section \ref{sec:numex} present preliminary numerical results. All proofs are relegated to the appendix.

\section{Problem statement}\label{sec:prb}
Throughout the paper, we shall assume that $K\subset E$ is a closed convex cone in Euclidean space $E$; we loose nothing by assuming that $K$ linearly spans $E$.
We assume, further, that $\|\cdot\|$ is a norm on $E$, and $f:K\to\bR$ is a convex function with Lipschitz continuous gradient, so that
$$\|f'(x)-f'(y)\|_*\leq L_f \|x-y\|\,\,\forall x,y\in K,$$
where $\|\cdot\|_*$ denotes the norm dual to $\|\cdot\|$, whence
\begin{equation}\label{suchthatinitial}
\forall x,y\in K: f(y)\leq f(x)+\langle f'(x),y-x\rangle +{L_f\over 2}\|y-x\|^2.
\end{equation}
We consider two kinds of problems, detailed below.

\paragraph{Norm-minimization.} Such problems correspond to
\begin{equation}
\label{prob1}
\rho_*=\min_{x} \left\{\|x\|:\,x\in K,\,f(x)\le 0\right\}.
\end{equation}
To tackle~\rf{prob1}, we consider the following parametric family of problems
\begin{equation}
\label{optrho}
\Opt(\rho)=\min\{f(x):\;\|x\|\le \rho, \,x\in K\} \, .
\end{equation}
Note that whenever (\ref{prob1}) is feasible, which we assume from now on, we have
\begin{equation}
\label{probopt}
\rho_*=\min\{\rho\geq0:\; \Opt(\rho)\le 0\},
\end{equation}
and both problems (\ref{prob1}), (\ref{probopt}) can be solved.

Given a tolerance $\epsilon>0$, we want to find an $\epsilon$-solution to the problem, that is, a pair $\rho_\epsilon$, $x_\epsilon\in K$ such that
\begin{equation}
\label{goalgbspp}
\rho_\epsilon\leq\rho_*\;\mbox{and}\;x_\epsilon\in X_{\rho_\epsilon}\;\mbox{such that}\;f(x_\epsilon)\le \epsilon,
\end{equation}
where  $X_{\rho}:=\{x\in E:\;x\in K,\;\|x\|\le \rho\}$. Getting back to the problem of interest (\ref{prob1}), $x_\epsilon$ is then ``super-optimal'' and $\epsilon$-feasible:
\[
\|x_\epsilon\|\le \rho_\epsilon\le \rho_*,\;\;\;f(x_\epsilon)\le \epsilon.
\]
\paragraph{Penalized norm minimization.} These problems write as
\begin{equation}
\Opt=\min_{x} \left\{f(x)+\kappa\|x\|:\;x\in K\right\} \, .
\label{prob21}
\end{equation}
A equivalent formulation is
\begin{equation}
\Opt=\min_{x,r} \left\{F([x;r])=\kappa r+f(x):x\in K,\|x\|\leq r\right\}.
\label{prob2}
\end{equation}
We shall refer to~\rf{prob2} as the problem of {\em composite optimization}
(CO).
Given a tolerance $\epsilon>0$, our goal is to find an $\epsilon$-solution to
(\ref{prob2}), defined as a feasible solution $(x_\epsilon, r_\epsilon)$ to the
problem satisfying $F([x_\epsilon;r_\epsilon])-\Opt\leq\epsilon$. Note that in
this case $x_\epsilon$ is an $\epsilon$-solution, in the similar sense, to
\rf{prob21}.
\paragraph{Special case.} In many applications where Problem (\ref{prob1}) arise,
(\ref{prob21}) the function $f$ enjoys a special structure:
$$
f(x) = \phi(\cA x-b),
$$
where $x\mapsto \cA x-b$ is an affine mapping from $E$ to $\bR^m$, and
$\phi(\cdot):\bR^m\to\bR$ is a convex function with Lipschitz continuous
gradient; we shall refer to this situation as to {\sl special case}. In such case,
the quantity $L_f$ can be bounded as follows.
Let
$\pi(\cdot)$ be some norm on $\bR^m$, $\pi_*(\cdot)$ be the conjugate norm, and
$\|\cA\|_{\|\cdot\|,\pi}$ be the norm of the linear mapping $x\mapsto \cA x$
{induced by the norms $\|\cdot\|$, $\pi(\cdot)$ on the argument and the image spaces:}
$$
\|\cA\|_{\|\cdot\|,\pi(\cdot)}=\max_{x\in E}\{\pi(\cA x):\|x\|\leq1\}.
$$
Let also $L_{\pi(\cdot)}[\phi]$ be the Lipschitz constant of the gradient of
$\phi$ induced by the norm $\pi(\cdot)$, so that
$$
\pi_*(\phi'(y)-\phi'(y')) \leq L_{\pi(\cdot)}[\phi]\pi(y-y')\,\,\forall y,y'\in
F.
$$
Then, one can take as $L_f$ the quantity
\begin{equation}\label{immediatelyseen}
L_f=L_{\pi(\cdot)}[\phi] \|\cA\|_{\|\cdot\|,\pi(\cdot)}^2.
\end{equation}
\par\noindent
{\bf Example 1: quadratic fit.} In many applications, we are interested in
$\|\cdot\|_2$-discrepancy between $\cA x$ and $b$; the related choice of
$\phi(\cdot)$ is $\phi(y)={1\over 2}y^Ty$. Specifying $\pi(\cdot)$ as
$\|\cdot\|_2$, we get $L_{\|\cdot\|_2}[\phi]=1$.\par\noindent
{\bf Example 2: smoothed $\ell_\infty$ fit.} When interested in
$\|\cdot\|_\infty$ discrepancy between $\cA x$ and $b$, we can use as $\phi$ the
function $\phi(y)={1\over 2} \|y\|_\beta^2$, where $\beta\in[2,\infty)$.
Taking $\pi(\cdot)$ as $\|\cdot\|_\infty$, we get
$$L_{\|\cdot\|_\infty}[\phi]\leq (\beta-1)m^{2/\beta}.$$ Note that
$$
{1\over 2}\|y\|_\infty^2\leq \phi(y)\leq {m^{2/\beta}\over2} \|y\|_\infty^2,
$$
so that for $\beta=O(1)\ln(m)$ and $m$ large enough (specifically, such that
$\beta\geq2$), $\phi(y)$ is within absolute constant factor of ${1\over
2}\|y\|_\infty^2$. The latter situation can
be interpreted as $\phi$ behaving as ${1\over2}\|\cdot\|_\infty^2$). At the same time, with
$\beta=O(1)\ln(m)$, $L_{\|\cdot\|_\infty}[\phi]\leq O(1)\ln(m)$ grows with $m$
logarithmically.
\par
Another widely used choice of $ \phi(\cdot)$ for this type of discrepancy is
``logistic'' function
\[
\phi(y)={1\over \beta}\ln\left(\sum_{i=1}^m \left[e^{\beta y_i}+e^{-\beta
y_i}\right]\right).
\]
 For $\pi(\cdot)=\|\cdot\|_{\infty}$ we easily compute
 $L_{\|\cdot\|_\infty}[\phi]\leq \beta$ and $\|y\|_\infty\le \phi(y)\le
\|y\|_\infty+ \ln(2n)/\beta$.
 \par
 Note that in some applications we are interested in ``one-sided'' discrepancies
quantifying the magnitude of the vector $[\cA x - b]_+=[\max[0,(\cA
x-b)_1];...;\max[0,(\cA x -b)_m]]$ rather than the the magnitude of the vector
$\cA x -b$ itself. Here, instead of using $\phi(y)=\half \|y\|^2_\beta$ in the
context of examples  1 and 2, one can use the functions $\phi_+(y)=\phi([y]_+)$.
In this case the bounds on $L_{\pi(\cdot)}[\phi_+]$ are exactly the same as the
above bounds on $L_{\pi(\cdot)}[\phi]$. The obvious substitute for
the two-sided logistic function is its ``one-sided version:'' $\phi_+(y)={1\over
\beta}\ln\left(\sum_{i=1}^m \left[e^{\beta y_i}+1\right]\right)$ which obeys the
same bound for $L_{\pi(\cdot)}[\phi_+]$ as its two-sided analogue.

\paragraph{First-order and Linear Optimization oracles.}
 We assume that $f$  is represented by a {\em first-order oracle} -- a routine
which, given on input a point $x\in K$, returns the value $f(x)$ and the
gradient $f'(x)$ of $f$ at $x$. As about $K$ and $\|\cdot\|$, we assume that
they are given by a {\sl Linear Optimization (LO) oracle} which, given on input a linear
form $\langle \eta,\cdot\rangle$ on $E$, returns a minimizer $x[\eta]$ of this
linear form on the set $\{x\in K:\|x\|\leq1\}$. We assume w.l.o.g. that for
every $\eta$, $x[\eta]$ is either zero, or is a vector of the $\|\cdot\|$-norm
equal to 1. To ensure this property, it suffices to compute $\langle
\eta,x[\eta]\rangle$ for $x[\eta]$ given by the oracle; if this inner product is
0, we can reset $x[\eta]=0$, otherwise $\|x[\eta]\|$ is automatically equal to
1.
 \par
 Note that {an LO oracle} for $K$ and $\|\cdot\|$ allows to find a minimizer
of a linear form of $z=[x;r]\in E^+:=E\times \bR$ on a set of the form
$K^+[\rho]=\{[x;r]\in E^+: x\in K, \|x\|\leq r\leq\rho\}$ due to the following observation:
\begin{lemma}\label{lem:evident} Let $\rho\geq0$ and
$\eta^+=[\eta;\sigma]\in E^+$. Consider the linear form $\ell(z)=\langle
\eta^+,z\rangle$ of $z=[x;r]\in E^+$, and let
$$
z^+=\left\{\begin{array}{ll}\rho[x[\eta];1]&,\langle
\eta^+,[x[\eta];1]\rangle\leq 0,\\
0&,\hbox{otherwise}
\\
\end{array}\right..
$$
Then $z^+$ is a minimizer of $\ell(z)$ over $z\in K^+[\rho]$, When $\sigma=0$,
one has $z^+=\rho[x[\eta];1]$.
\end{lemma}
Indeed, let $z_*=[x^*;r_*]$ be a minimizer of $\ell(\cdot)$ over $K^+[\rho]$. Since $\|x^*\|\leq r_*$ due to $[x^*;r_*]\in K^+[\rho]$, we have
$z^*:=r^*[x[\eta];1]\in K^+[\rho]$ due to $\|x[\eta]\|\leq1$, and $\ell*(z^*)\leq\ell(z_*)$ due to the definition of $x[\eta]$. We conclude that
any minimizer of $\ell(\cdot)$ over the segment $\{s[x[\eta];1]:0\leq s\leq\rho\}$ is also a minimizer of $\ell(\cdot)$ over $K^+[\rho]$. It remains to note that the vector indicated in Lemma clearly is a minimizer of $\ell(\cdot)$ on the above segment. \qed

 \section{Conditional Gradient algorithm}\label{sec:po}
%
In this section, we present an overview of the properties of the standard Conditional Gradient algorithm,
and highlight some memory-based extensions. These properties are not new. However,
since they are key for the design of our proposed algorithms in the next sections,
we present them for further reference.

\subsection{Conditional gradient algorithm}\label{subs:CG}
Let $E$ be a Euclidean space and $X$ be a closed and bounded convex set in $E$
which linearly spans $E$. Assume that $X$ is given by a LO oracle -- a
routine which, given on input $\eta\in E$, returns an optimal solution
$x_X[\eta]$ to the optimization problem
\begin{equation*}
\min_{x\in X} \: \langle \eta,x\rangle
\end{equation*}
(cf. Section \ref{sec:prb}).
Let $f$ be a convex differentiable function on $X$ with Lipschitz continuous
gradient $f'(x)$, so that
\begin{equation}
\label{suchthat}
\forall x,y\in X: f(y)\leq f(x)+\langle f'(x),y-x\rangle +{\half}L\|y-x\|_X^2,
\end{equation}
where $\|\cdot\|_X$ is the norm on $E$ with the unit ball $X-X$. We intend to
solve the problem
\begin{equation}\label{problemf}
f_*=\min_{x\in X} f(x).
\end{equation}
A {\sl generic {\CG} algorithm} is a recurrence which builds iterates $x_t \in X$, $t=1,2,...$, in such a way that
\begin{equation}\label{newneweq1}
f(x_{t+1})\leq f(\widetilde{x}_{t+1}),
\end{equation}
where
\begin{equation}\label{newneweq2}
\begin{array}{rcl}
\widetilde{x}_{t+1}&=&x_t+\gamma_t[x_t^+-x_t],\;\;\hbox{\ where\ }
x_t^+=x_X[f'(x_t)] \hbox{\ and\ }\gamma_t={2\over t+1}
.\\
\end{array}
\end{equation}
Basic implementations of a generic {\CG} algorithm are given by
\begin{equation}\label{eqreq}
\begin{array}{lrcl}
(a)&x_{t+1}&=&x_t+\gamma_t[x_t^+-x_t],\;\;\gamma_t={2\over t+1},\\
(b)&x_{t+1}&\in&\Argmin_{x\in D_t} f(x),\;\; D_t=[x_t,x_t^+];
\end{array}
\end{equation}
in the sequel, we refer to them as {\CG}a and {\CG}b, respectively. As a byproduct of running generic {\CG}, after $t$
steps we have at our disposal the quantities
\begin{equation}\label{eq:lowbnd}
f_{*,k}=\min\limits_{x\in X}\left[f(x_k)+\langle
f'(x_k),x-x_k\rangle\right]=f(x_k)-\langle
f'(x_k),x_k-x_X[f'(x_k)]\rangle,\,\,1\leq k\leq t,
\end{equation}
which, by convexity of $f$, are lower bounds on $f_*$. Consequently, at the end
of step $t$ we have at our disposal a lower bound
\begin{equation}\label{eq:lowbnd1}
f_*^t:=\max_{1\leq k\leq t} f_{*,k}\leq f_*,\,t=1,2,...
\end{equation}
on $f_*$.
\par Finally, we define the approximate solution $\bar{x}_t$ found in course of
$t=1,2,...$ steps  as the best -- with the smallest value of $f$ -- of the
points
$x_1,...,x_t$. Note that $\bar{x}_t\in X$.
\par
{The following statement summarizes the well known properties of {\CG} (to make the
presentation self-contained, we provide in Appendix the proof).}
\begin{theorem}\label{propredgrad} For a generic {\CG} algorithm, in particular, for both  {\CG}a, {\CG}b, we have
\begin{equation}\label{then}
f(\bar{x}_t)-f_*\leq f(x_{t})-f_*\leq{2L\over t+1},\quad t\geq 2 ;
\end{equation}
and
\begin{equation}
\label{certif1}
f(\bar{x}_t)-f_*^t\leq {4.5 L\over t-2},{\quad} t\geq 5 .
\end{equation}
\end{theorem}
Some remarks regarding the conditional algorithm are in order.

{\bf Certifying quality of approximate solutions.} An attractive property of {\CG}
is the presence of online lower bound $f_*^t$ on $f_*$ which certifies the
theoretical rate of convergence of the algorithm, see (\ref{certif1}). This accuracy certificate, first established in \cite{Jaggi:2013},
also provides a valuable stopping criterion when running the algorithm in practice.
\par
{\bf {\CG} algorithm with memory.} When computing the next search point $x_{t+1}$
the simplest {\CG} algorithm {\CG}a only uses the latest answer
$x_t^+=x_X[f'(x_t)]$ of the {LO oracle}. Meanwhile, algorithm {\CG}b can be
modified to make use of information supplied by previous oracle calls; we refer
to this modification as {\CG} with memory ({\CG}M).\footnote{Note that in the context
of ``classical'' Frank-Wolfe algorithm -- minimization of a smooth function over
a polyhedral set -- such modification is referred to as {\em Restricted
Simplicial Decomposition} \cite{Hol74,HeLaVe87,VeHe93}}. Assume that we have
already carried out $t-1$ steps of the algorithm and have at our disposal
current iterate $x_t\in X$ (with $x_1$ selected as an arbitrary point of $X$)
along with previous iterates $x_\tau$, $\tau<t$ and the vectors $f'(x_\tau)$,
$x_\tau^+=x_X[f'(x_\tau)]$. At the step, we compute $f'(x_t)$ and
$x^+_t=x_X[f'(x_t)]$. Thus, at this point in time we have at our disposal $2t$
points $x_\tau,x^+_\tau$, $1\leq \tau\leq t$, which belong to $X$. Let $X_t$ be
subset of these points, with the only restriction that the points $x_t$, $x^+_t$
are selected, and let us define the next iterate $x_{t+1}$ as
\begin{equation}\label{bundlex}
x_{t+1}\in\Argmin\limits_{x\in \Conv(X_t)} f(x),
\end{equation}
that is,
\begin{equation}\label{bundleaux}
x_{t+1}=\sum_{x\in X_t}\lambda^t_x
x,\,\,\lambda^t\in\Argmin\limits_{\lambda^t=\{\lambda_x\}_{x\in X_t}}
\left\{f\left({\sum}_{x\in X_t} \lambda_x x\right):\;\lambda\geq0,\;\sum_{x\in
X_t}\lambda_x=1\right\}.
\end{equation}
Clearly, it is again a generic {\CG} algorithm, so that conclusions in Theorem \ref{propredgrad} are fully applicable to {\CG}M.
Note that {\CG}b per se is nothing but
{\CG}M  with $X_t=\{x_t,x_t^+\}$ {and $M=2$} for all $t$.\par\noindent
{\bf {\CG}M: implementation issues.} Assume that the cardinalities of the sets
$X_t$ in {\CG}M are bounded by some {$M \geq 2$}.   In this case, implementation of the method requires solving
at every step an auxiliary problem (\ref{bundleaux}) of minimizing  over the
standard simplex of dimension $\leq M-1$ a smooth convex function given by a
first-order oracle {induced by the first-oracle for $f$}.
When $M$ is a once for ever fixed small integer, the arithmetic cost of solving
this problem within machine accuracy by, say, the Ellipsoid algorithm is
dominated by the arithmetic cost of just $O(1)$ calls to the first-order oracle
for $f$. Thus, {\CG}M with small $M$ can be considered as
implementable\footnote{Assuming possibility to solve (\ref{bundleaux}) exactly,
while being idealization, is basically as ``tolerable'' as the standard in
continuous optimization assumption that one can use exact real arithmetic or
compute exactly eigenvalues/eigenvectors of symmetric matrices. The outlined
``real life'' considerations can be replaced with rigorous error analysis which
shows that in order to maintain the efficiency estimates from Theorem
\ref{propredgrad}, it suffices to solve $t$-th auxiliary problem within properly
selected positive inaccuracy, and this can be achieved in $O(\ln(t))$
computations of $f$ and $f'$.}.

Note that in the special case (Section \ref{sec:prb}), where $f(x)=\phi(Ax-b)$,
assuming $\phi(\cdot)$ and $\phi'(\cdot)$ easy to compute, as is the case in
most of the applications, the first-order oracle for the auxiliary problems
arising in {\CG}M becomes cheap (cf.  \cite{ZibulNarkiss}). Indeed, in this case
(\ref{bundleaux}) reads
$$
\min_{\lambda^t} \left\{g_t(\lambda^t):=\phi\left({\sum}_{x\in
X_t}\lambda^t_xAx-b\right): \begin{array}{l}
\lambda^t=\{\lambda^t_x\}_{x\in X_t}\geq0,\\
{\sum}_{x\in X_t}\lambda^t_x=1\\
\end{array}\right\}.
$$
It follows that all we need to get a computationally cheap access to the
first-order information on $g_t(\lambda^t)$ {\sl for all} values of $\lambda^t$
is to have at our disposal the matrix-vector products $Ax$, $x\in X_t$. With our
construction of $X_t$, the only two ``new'' elements in $X_t$ (those which were
not available at preceding iterations) are $x_t$ and $x^+_t$, so that the only
two new matrix-vector products we need to compute at iteration $t$ are $Ax_t$
(which usually is a byproduct of computing $f'(x_t)$) and $Ax_t^+$. Thus, we can
say that the ``computational overhead,'' as compared to computing $f'(x_t)$ and
$x_t^+=x_X[f'(x_t)]$, needed to get easy access to the first-order information
on $g_t(\cdot)$ reduces to computing the single matrix-vector product $Ax_t^+$.

\section{Conditional gradient algorithm for parametric optimization}\label{subsec:posec}
{In this section, we describe a multi-stage algorithm to solve the parametric optimization problem (\ref{optrho}), (\ref{probopt}), using the conditional algorithm to solve inner sub-problems.}
(\ref{optrho}), (\ref{probopt}). The idea, originating from \cite{Lem:Nem:Nes:1995} (see also \cite{Nesterov:2004,NemJud1,Nes:Nem:2013}), is to use a Newton-type method for approximating from below the positive root $\rho_*$ of $\Opt(\rho)$, with (inexact) first-order information on $\Opt(\cdot)$ yielded by approximate solving the optimization problems defining $\Opt(\cdot)$; the difference with the outlined references is that now we solve these problems with the  {\CG} algorithm.
\par
Our algorithm works stagewise.
At the beginning of stage $s=1,2,...$, we have at hand a lower bound $\rho_s$ on $\rho_*$, with $\rho_1$  defined as follows:
\begin{quote}
We compute $f(0)$, $f'(0)$ and $x[f'(0)]$. If $f(0)\leq\epsilon$ or $x[f'(0)]=0$, we are done --- the pair ($\rho=0$, $x=0$) is an  $\epsilon$-solution to (\ref{probopt}) in the first case, and is an optimal solution to the problem in the second case (since in the latter case $0$ is a minimizer of $f$ on $K$, and (\ref{probopt}) is feasible). Assume from now on that the above options do not take place (``nontrivial case''), and let
$$
d=-\langle f'(0),x[f'(0)]\rangle.
$$
Due to the origin of $x[\cdot]$, $d$ is positive, and $f(x)\geq f(0)+\langle f'(0),x\rangle \geq f(0)-d\|x\|$ for all $x\in K$, which implies that $
\rho_*\geq \rho_1:={f(0)\over d}>0.$
\end{quote}
\par  At stage $s$ we apply a generic {\CG} algorithm (e.g., {\CG}a,{\CG}b, or {\CG}M; in the sequel, we refer to the algorithm we use as to {\CG}) to the auxiliary problem
\begin{equation}\label{eqauxpcg}
\Opt(\rho_s)=\min_x\{f(x):\,x\in K[\rho_s]\},\quad K[\rho]=\{x\in K:\|x\|\leq\rho\},
\end{equation}
Note that the {LO oracle} for $K$, $\|\cdot\|$ induces {an LO oracle} for $K[\rho]$; specifically, for every $\eta\in E$, the point $x_\rho[\eta]:=\rho x[\eta]$ is a minimizer of the linear form $\langle\eta,x\rangle$ over $x\in K[\rho]$, see Lemma \ref{lem:evident}. $x_\rho[\cdot]$ is exactly the {LO oracle} utilized by {\CG} as applied to (\ref{eqauxpcg}).

As explained above, after $t$ steps of {\CG} as applied to (\ref{eqauxpcg}), the iterates being $x_\tau\in K[\rho_s]$, $1\leq\tau\leq t$ \footnote{The iterates $x_t$, same as other indexed by $t$ quantities participating in the description of the algorithm, in fact depend on both $t$ and the stage number $s$. To avoid cumbersome notation when speaking about a particular stage, we suppress $s$ in the notation.},  we have at our disposal current approximate solution $\bar{x}_t\in\{x_1,...,x_t\}$ such that $f(\bar{x}_t)=\min_{1\leq\tau\leq t} f(x_\tau)$ along with a lower bound $f_*^t$ on $\Opt(\rho_s)$. Our policy is as follows.
\begin{enumerate}
\item When $f(\bar{x}_t)\leq\epsilon$, we terminate the solution process and output $\bar{\rho}=\rho_s$ and $\bar{x}=\bar{x}_t$;
\item When the above option is not met and $f_*^t<{3\over 4}f(\bar{x}_t)$, we specify $x_{t+1}$ according to the description of {\CG} and pass to step $t+1$ of stage $s$;
\item Finally, when neither one of the above options takes place, we terminate stage $s$ and pass to stage $s+1$, specifying $\rho_{s+1}$ as follows:\\
 We are in the situation $f(\bar{x}_t)>\epsilon$ and $f_*^t\geq {3\over 4}f(\bar{x}_t)$. Now, for $k\leq t$ the quantities $f(x_k)$, $f'(x_k)$ and $x[f'(x_k)]$ define affine function of $\rho\geq0$
 $$
 \ell_k(\rho)=f(x_k)+\langle f'(x_k),x-\rho x[f'(x_k)]\rangle.
 $$
By Lemma \ref{lem:evident} we have for every $\rho\geq0$
$$
\ell_k(\rho)=\min_{x\in K[\rho]}\left[f(x_k)+\langle f'(x_k),x-x_k\rangle\right]\leq \min_{x\in K[\rho]} f(x)=\Opt(\rho),
$$
where the inequality is due to the convexity of $f$. Thus, $\ell_k(\rho)$ is an affine in $\rho\geq0$ lower bound on $\Opt(\rho)$, and we lose nothing by assuming that all these univariate affine functions are memorized when running {\CG} on (\ref{eqauxpcg}). Note that by construction of the lower bound $f_*^t$ (see (\ref{eq:lowbnd}), (\ref{eq:lowbnd1}) and take into account that we are in the case of $X=K[\rho_s]$, $x_X[\eta]=\rho_s x[\eta]$) we have
$$f_*^t=\ell^t(\rho_s),\,\,\ell^t(\rho)=\max_{1\leq k\leq t}\ell_k(\rho).$$
Note that $\ell^t(\rho)$ is a lower bound on $\Opt(\rho)$, so that $\ell^t(\rho)\leq0$ for $\rho\geq\rho_*$, while $\ell^t(\rho_s)=f_*^t$ is positive. It follows that
 $$
 r^t:=\min\left\{\rho:\ell^t(\rho)\leq0\right\}
$$
is well defined and satisfies $\rho_s<r^t\leq\rho_*$. We compute $r^t$ (which is easy) and pass to stage $s+1$, setting $\rho_{s+1}=r^t$ and selecting, as the first iterate of the new stage, any point known to belong to $K[\rho]$ (e.g., the origin, or $\bar{x}_t$). The first iterate of the first stage is $0$.
\end{enumerate}
The description of the algorithm is complete.

The complexity properties of the algorithm are given by the following proposition.
\begin{theorem}\label{TheMain2}
When solving a PO problem \rf{optrho}, \rf{probopt} by the outlined algorithm,
\par
{\rm (i)} the algorithm terminates with an $\epsilon$-solution,
as defined in Section \ref{sec:prb} (cf. \rf{goalgbspp});
\par
{\rm (ii)} The number $N_s$ of steps at every stage $s$ of the method admits the bound
\[
N_s\le \max\left[6,{72\rho_*^2L_f\over\epsilon}+3\right].
\]
\par
{\rm (iii)} The number of stages before termination does not exceed the quantity
$$\max\left[1.2\ln\left({f(0)+\half L_f \rho_*^2\over\epsilon^2}\right)+2.4,3\right].$$

\end{theorem}

 \section{Conditional Gradient algorithm for Composite Optimization}\label{sec:cocg}
{In this section, we present a modification of the \CG~algorithm capable to solve composite minimization problem (\ref{prob2}). We assume in the sequel that $\|\cdot\|,K$ are represented by an LO oracle for the set $\{x\in K:\|x\|\leq1\}$, and $f$ is given by a first order oracle.}
In order to apply {\CG} to the composite optimization problem \rf{prob2}, we make the assumption as follows:
\begin{quote} {\bf Assumption A:}
There exists $D<\infty$ such that $\kappa r+f(x)\leq f(0)$ together with $\|x\|\leq r$, $x\in K$, imply that $r\leq D$.
\end{quote}
We define $D_*$ as the minimal value of $D$ satisfying Assumption A, and assume that {\sl we have at our disposal a finite upper bound $D^+$ on $D_*$.} An important property of the algorithm we are about to develop is that its efficiency estimate depends on the induced by problem's data quantity $D_*$, and is independent of our a priori upper bound $D^+$ on this quantity, see Theorem \ref{pro:cgco} below.
\paragraph{The algorithm.} We are about to present an algorithm for solving \rf{prob2}.
Let $E^+=E\times\bR$, and $K^+=\{[x;r]:\;x\in K,\,\|x\|\le r\}$.
From now on, for a point $z=[x;r]\in E^+$ we set   $x(z)=x$ and $r(z)=r$.
Given $z=[x;r]\in K^+$, let us consider the segment
\[
\Delta(z)=\{\rho[x[f'(x)];1]:\;0\leq \rho\leq D^+\}.
\]
and
the linear form
\[
{ \zeta=[\xi;\tau]\to \langle f'(x),\,\xi\rangle +\kappa \tau= \langle F'(z),\,\zeta\rangle}
 \]
Observe that by Lemma \ref{lem:evident}, for every $0\leq \rho\leq D^+$, the minimum of this form on $K^+[\rho]=\{[x;r]\in E^+, x\in K, \|x\|\leq r\leq\rho\}$  is attained at a point of $\Delta(z)$ (either at
$[\rho x[f'(x)];\,\rho]$ or at the origin).
A generic Conditional Gradient algorithm for composite optimization (CO{\CG}) is a recurrence which builds the points $z_t=[x_t;r_t]\in K^+$, $t=1,2,...$, in such a way that
\begin{equation}
 \label{cgco}
z_1=0;\,F(z_{t+1})\leq \min_z\{F(z):\;z\in \Conv\left(\Delta(z_t)\cup\{z_t\}\right)\},\;\;t=1,2,...
\end{equation}
Let $z_*=[x_*;r_*]$ be an optimal solution to \rf{prob2} (which under Assumption A clearly exists), and let $F_*=F(z_*)$ (i.e., $F_*$ is nothing but $\Opt$, see (\ref{prob21})).
\begin{theorem}
\label{pro:cgco}
A generic CO{\CG} algorithm  \rf{cgco} maintains the inclusions $z_t\in K^+$ and is a descent algorithm: $F(z_{t+1})\leq F(z_t)$ for all $t$. Besides this, we have
\begin{equation}\label{eq:cocgrate}
F(z_t)-F_*\le {8L_f D_*^2\over t+14},\;t=2,3, ...
\end{equation}
\end{theorem}

\paragraph{CO{\CG} with memory.}
The simplest implementation of a generic CO{\CG} algorithm is given by the recurrence
\begin{equation}
\label{cgco_simple}
z_1=0;\,z_{t+1}\equiv [x_{t+1};r_{t+1}]\in \Argmin_z\{F(z):\;z\in \Conv\left(\Delta(z_t)\cup\{z_t\}\right)\},\;\;t=1,2,...\;.
\end{equation}
Denoting $\widehat{z}_\tau :=D^+[x[f'(x_\tau)];1]$, the recurrence can be written
\begin{equation}
\label{aux}
\begin{array}{l}
z_{t+1}=\lambda_t \widehat{z}_{t} +\mu_tz_t,\hbox{\ where}\\
(\lambda_t,\mu_t)\in \Argmin\limits_{\lambda,\mu}\bigg\{F(\lambda \widehat{z}_{t} +\mu z_t):\;\lambda+\mu\leq1,\;\lambda\geq0,\mu\geq0\bigg\}.
\end{array}
\end{equation}
 As for the {\CG} algorithm in section~\ref{sec:po}, the recurrence \rf{cgco_simple} admits a version with memory CO{\CG}M still obeying (\ref{cgco})  and thus sartisfying the conclusion of Theorem \ref{pro:cgco}. Specifically, assume that we already have built $t$ iterates $z_\tau=[x_\tau;r_\tau]\in K^+$, $1\leq\tau\leq t$, with $z_1=0$, along with the gradients $f'(x_\tau)$ and the points $x[f'(x_\tau)]$. Then we have at our disposal a number of points from $K^+$, namely, the iterates $z_\tau$, $\tau\leq t$, and the points $\widehat{z}_\tau=D^+[x[f'(x_\tau)];1]$. Let us select a subset $Z_t$ of the set $\{z_\tau,\widehat{z}_\tau,1\leq \tau\leq t\}$, with the only restriction that $Z_t$ contains the points $z_t,\widehat{z}_t$, and set \begin{equation} \label{bundle3} z_{t+1}\in\Argmin_{z\in \C_t} F(z), \quad  \C_t=\Conv\{\{0\}\cup Z_t\}\}.
\end{equation}
Since $z_t,\widehat{z}_t\in Z_t$, we have $\Conv\left(\Delta(z_t)\cup\{z_t\}\right)\}\subset\C_t$, whence the procedure we have outlined is an implementation of  generic CO{\CG} algorithm. Note that the basic CO{\CG} algorithm is the particular case of the CO{\CG}M corresponding to the case where $Z_t=\{z_t,\widehat{z}_t\}$ for all $t$. The discussion of implementability of  {\CG}M in section \ref{sec:po} fully applies to CO{\CG}M.
\par
  Let us outline several options which can be  implemented in CO{\CG}M; while preserving the theoretical efficiency estimates stated in Theorem \ref{pro:cgco}
they can improve the practical performance of the algorithm. For the sake of definiteness, let us focus on the case of quadratic $f$: $f(x)=\|\cA x- b\|_2^2$, with $\Ker \cA =\{0\}$; extensions to a more general case are straightforward.
\begin{enumerate}
\item[{\bf A.}] We lose nothing (and potentially gain) when extending $\C_t$ in
(\ref{bundle3}) to  the conic hull
$$
\C_t^+=\{w=\sum_{\zeta\in Z_t}\lambda_\zeta \zeta:\;\lambda_\zeta\geq0,\;\zeta\in Z_t\} $$ of $Z_t$. When $K=E$, we can go further and replace (\ref{bundle3}) with \begin{equation}\label{bundleprime}
z_{t+1}\in\Argmin_{z=[x;r],\lambda}\left\{f(x)+\kappa
r:\;x=\sum_{\zeta=[\eta;\rho]\in Z_t}\lambda_\zeta \eta,\;r\geq \sum_{\zeta=[\eta;\rho]\in Z_t}|\lambda_\zeta|\rho\right\}.
\end{equation}
Note that the preceding ``conic case'' is obtained from
(\ref{bundleprime}) by adding to the constraints of the right hand side problem the inequalities $\lambda_\zeta\geq0,\zeta\in Z_t$. Finally, when  $\|\cdot\|$ is easy to compute, we can improve (\ref{bundleprime}) to
\begin{equation}\label{bundledoubleprime}
\begin{array}{c}
z_{t+1}=\left[\sum_{\zeta=[\eta;\rho]\in Z_t}\lambda^*_\zeta\eta;\left\|\sum_{\zeta=[\eta;\rho]\in Z_t}\lambda_\zeta^*\eta\right\|\right],\\
\lambda^*\in\Argmin_{\{\lambda_\zeta,\zeta\in Z_t\}}\left\{f\left(\sum_{\zeta=[\eta;\rho]\in Z_t}\lambda_\zeta\eta\right)+\kappa
\sum_{\zeta=[\eta;\rho]
\in Z_t} |\lambda_\zeta|\rho\right\}\\
\end{array}
\end{equation}
(the definition of $\lambda^*$ assumes that $K=E$, otherwise the constraints of the problem specifying $\lambda^*$ should be augmented by the inequalities $\lambda_\zeta\geq0,\zeta\in Z_t$).
\item[{\bf B.}] In the case of quadratic $f$ and moderate cardinality of $Z_t$, optimization problems arising in  (\ref{bundleprime}) (with or without added constraints $\lambda_\zeta\geq0$) are explicitly given low-dimensional ``nearly quadratic'' convex problems which can be solved to high accuracy ``in no time'' by interior point solvers. With this in mind, we could solve these problems for the given value of the penalty parameter $\kappa$ {\sl and also for several other values of the parameter}.
Thus, at every iteration we get feasible approximate solution to several  instances of (\ref{prob21}) for different values of the penalty parameter. Assume that we keep in memory, for every value of the penalty parameter in question, the best, in terms of the respective objective, of the related approximate solutions found so far. Then upon termination we will have at our disposal, along with the feasible approximate solution associated with the given value of the penalty parameter, provably obeying the efficiency estimates of Theorem \ref{pro:cgco}, a set of feasible approximate solutions to the instances of (\ref{prob21}) corresponding to other values of the penalty.
\item[{\bf C.}] In the above description, $Z_t$ was assumed to be a subset of the set $Z^t=\{z_\tau=[x_\tau;r_\tau],\widehat{z}_\tau,\,1\leq\tau\leq t\}$ containing $z_t$ and $\widehat{z}_t$. Under the latter restriction, we lose nothing when allowing for $Z_t$ to contain points from $K^+\backslash Z^t$ as well. For instance, when $K=E$ and $\|\cdot\|$ is easy to compute, we can add to $Z_t$ the point $z_t^\prime=[f'(x_t);\|f'(x_t)\|]$. Assume, e.g., that we fix in advance the cardinality $M\geq3$ of $Z_t$ and define $Z_t$ as follows: to get $Z_t$ from $Z_{t-1}$, we eliminate from the latter set several (the less, the
better) points to get a set of cardinality $\leq M-3$, and then add to the resulting set the points $z_t$, $\widehat{z}_t$ and $z_t^\prime$.
Eliminating the points according to the rule ``first in -- first out,''
the projection of the feasible set of the optimization problem in
(\ref{bundledoubleprime})  onto the space of $x$-variables will be a linear subspace of $E$ containing, starting with step $t=M$, at least $\lfloor M/3\rfloor$ (here $\lfloor a \rfloor$ stands for the largest integer not larger than $a$) of  gradients of $f$ taken at the latest iterates, so that the method, modulo the influence of the penalty term, becomes a ``truncated'' version of the Conjugate Gradient algorithm for quadratic minimization. Due to nice convergence properties of Conjugate Gradient in the quadratic case, one can hope that a modification of this type will improve significantly the practical performance of CO{\CG}M.
\end{enumerate}

\section{Application examples}\label{sec:examples}
{In this section, we detail how the proposed conditional gradient algorithms apply to several examples. In particular,
we detail the corresponding {LO} oracles, and how one could implement these oracles efficiently. }
\subsection{Regularization by nuclear/trace norm}
The first example where the proposed algorithms seem to be more attractive than the proximal methods are large-scale problems (\ref{prob1}), (\ref{prob21}) on the space of $p\times q$ matrices  $E=\bR^{p\times q}$
{associated
with the nuclear norm $\Vert \sigma(x) \Vert_1$ of a matrix $x$, where $\sigma(x)=[\sigma_1(x);...;\sigma_{\min[p,q]}(x)]$}
is the
vector of singular values of a $p\times q$
matrix $x$. Problems of this type with $K=E$ arise in various versions of matrix
completion, where the goal is to recover a matrix $x$ from its noisy linear
image  $y=\cA x +\xi$, so that $f=\phi(\cA x-y)$, with some smooth and convex
discrepancy measure $\phi(\cdot)$, most notably, $\phi(z)={1\over 2}\|z\|_2^2$.
In this case, $\|\cdot\|$ minimization/penalization is aimed at getting a
recovery of low rank
(\cite{SrebroShraibman2005,candes-2008,care2009exact,
Goldfarb:2009,Jaggi:2010,parillo-2010,RechtRe2011,YaYa2011,Ma:Goldf:2011,
ShalevShwartzGoSh11}
 and references therein). Another series of applications relates to the case
when $E=\bS^p$ is the space of symmetric $p\times p$ matrices, and $K=\bS^p_+$
is the cone of positive semidefinite matrices, with $f$ and $\phi$ as above;
this setup corresponds to the situation when one wants to recover a covariance
(and thus positive semidefinite symmetric) matrix from experimental data.
Restricted from $\bR^{p\times p}$ onto $\bS^p$, the nuclear norm becomes the
trace norm $\|\lambda(x)\|_1$, where $\lambda(x)\in\bR^p$ is the vector of
eigenvalues of a symmetric $p\times p$
matrix $x$, and regularization by this norm is, as above,  aimed at building a
low rank recovery.

\par
With the nuclear (or trace) norm in the role of $\|\cdot\|$, all known proximal
algorithms require, at least in theory, computing at every iteration the
complete singular value decomposition of $p\times q$ matrix $x$ (resp., complete
eigenvalue decomposition of a symmetric $p\times p$ matrix  $x$), which for
large $p,q$ may become prohibitively time consuming. In contrast to this, with
$K=E$ and {$\|\cdot\|=\|\sigma(\cdot)\|_1$, LO oracle for $(K,\|\cdot\|=\|\sigma(\cdot)\|_1)$} only requires
computing the leading right singular vector $e$ of a $p\times q$ matrix $\eta$
(i.e., the leading eigenvector of $\eta^T\eta$): $x[\eta]=-\bar{f}\bar{e}^T$,
where $\bar{e}=e/\|e\|_2$ and $\bar{f}=\eta e/\|\eta e\|_2$ for nonzero $\eta$
and $\bar{f}=0$, $\bar{e}=0$ when $\eta=0$. Computing the leading singular
vector of a large matrix is, in most cases, much cheaper than computing the
complete eigenvalue decomposition of the matrix. Similarly, in the case of
$E=\bS^p$, $K=\bS^p_+$ and  the trace norm in the role of $\|\cdot\|$, {LO oracle} requires computing the leading  eigenvector $e$ of a matrix
$\eta\in \bS^p$: $x[-\eta]=\bar{e}\bar{e}^T$, where $\bar{e}=0$ when $e^T\eta
e\geq0$, and $\bar{e}=e/\|e\|_2$ otherwise. Here again, for a large symmetric
$p\times p$ matrix, the required computation usually is much easier than
computing the complete eigenvalue decomposition of such a matrix. As a result,
in the situations under consideration, algorithms based on the LO
oracle remain ``practically implementable'' in an essentially larger range of
problem sizes than proximal methods.
\par
An additional attractive property of the {\CG} algorithms we have described stems
from the fact that {\sl since in the situations in question the matrices
$x[\eta]$ are of rank 1, $t$-th approximate solution $x_t$ yielded by the {\CG}
algorithms for composite minimization from Section \ref{sec:cocg} is of rank at
most $t$. Similar statement holds true for $t$-th approximate solution $x_t$
built at a stage of a {\CG} algorithm for parametric optimization from Section
\ref{sec:po}, provided that the first iterate at every stage is {the zero matrix.}}\footnote{this property is an immediate corollary of the fact that in the
situation in question, by description of the algorithms $x_t$ is a convex
combination of $t$ points of the form $x[\cdot]$.}.
\subsection{Regularization by Total Variation}\label{sect:TV}
Given integer $n\geq2$, consider the linear space $M^n:=\bR^{n\times n}$.  We interpret elements $x$ of $M^n$ as {\sl images} -- real-valued functions $x(i,j)$ on the $n\times n$ grid $\Gamma_{n,n}=\{[i;j])\in\bZ^2: 0\leq i,j<n\}$. The (anisotropic) Total Variation (TV) of an image $x$ is the $\ell_1$-norm of its (discrete) gradient field $(\nabla_ix(\cdot),\nabla_jx(\cdot))$:
$$
\begin{array}{c}
\TV(x)=\|\nabla_i x\|_1+\|\nabla_j x\|_1,\\
\begin{array}{rcl}\nabla_ix(i,j)&=&x_(i+1,j)-x(i,j):\Gamma_{n-1,n}:=\{[i;j]\in\bZ^2:0\leq i<n-1,0\leq j<n\},\\
\nabla_jx(i,j)&=&x_(i,j+1)-x(i,j):\Gamma_{n,n-1}:=\{[i;j]\in\bZ^2:0\leq i<n,0\leq j<n-1\}\\
\end{array}
\end{array}
$$
Note that $\TV(\cdot)$ is a norm on the subspace $M^n_0$ of $M^n$ comprised of {\sl zero mean images} $x$ (those with $\sum_{i,j}x(i,j)=0$) and vanishes on the orthogonal complement to $M_0^n$, comprised of constant images. \par
Originating from the celebrated paper \cite{Osheretal} and extremely popular {\sl Total Variation-based image reconstruction} in its basic version recovers an image $x$ from its noisy observation $b=\cA x +\xi$ by solving
problems (\ref{prob1}) or (\ref{prob21}) with $K=E=M^n$, $f(x)=\phi(\cA x-b)$ and the seminorm $\TV(\cdot)$ in the role of $\|\cdot\|$. In the sequel, we focus on the versions of these problems where $K=E=M^n$ is replaced with $K=E=M^n_0$, thus bringing the $TV$-regularized problems into our framework. This restriction is basically harmless; for example, in the most popular case of $f(x)={1\over 2}\|\cA x-b\|_2^2$ reduction to the case of $x\in M^n_0$ is immediate -- it  suffices to replace $(\cA,b)$ with $(P\cA,Pb)$, where $P$ is the orthoprojector onto the orthogonal complement to the one-dimensional subspace spanned by $\cA{\mathbf{e}}$, where $\mathbf{e}$ is the all-ones image\footnote{When $f$ is more complicated, optimal adjustment of the mean $t$ of the image reduces by bisection in $t$ to solving small series of problems of the same structure as (\ref{prob1}), (\ref{prob21}) where the mean of the image $x$ is fixed and, consequently, the problems reduce to those with $x\in M^n_0$ by shifting $b$.}.
Now, large scale  problems (\ref{prob1}), (\ref{prob21}) with $K=E=M^n_0$ and $\TV(\cdot)$ in the role of $\|\cdot\|$ are difficult to solve by proximal algorithms. Indeed, in the situation in question a proximal algorithm would require at every iteration either minimizing function of the form $\TV(x)+\langle e,x\rangle +\omega(x)$ over  the entire $E$, or minimizing function of the form $\langle e,x\rangle +\omega(x)$ on a $\TV$-ball\footnote{which one of these two options takes place depends on the type of the algorithm.}, where $\omega(x)$ is albeit simple, but  {\sl nonlinear} convex function (e.g., $\|x\|_2^2$, or $\|\nabla_ix\|_2^2+\|\nabla_jx\|_2^2$). Auxiliary problems of this type seem to be difficult in the large scale case, especially taking into account that  when running  a proximal algorithm we need to solve at least tens, and more realistically -- hundreds of them\footnote{On a closest inspection, ``complex geometry'' of the $\TV$-norm stems from the fact that after parameterizing a zero mean image by its discrete gradient field and treating this field $(g=\nabla_i x,h=\nabla_jx)$ as our new design variable, the unit ball of the $\TV$-norm becomes the intersection of a simple set in the space of pairs $(g,h)\in F=\bR^{(n-1)\times n}\times \bR^{n\times(n-1)}$ (the $\ell_1$ ball $\Delta$ given by $\|g\|_1+\|h\|_1\leq 1$) with a linear subspace $P$ of $F$ comprised of {\sl potential} vector fields $(f,g)$ -- those which indeed are discrete gradient fields of images. Both dimension and codimension of $P$ are of order of $n^2$, which makes it difficult to minimize over $\Delta\cap P$ {\sl nonlinear}, even simple, convex functions, which is exactly what is needed in proximal methods.}.    In contrast to this, a LO oracle for the unit ball $\CTV=\{x\in M^n_0:\TV(x)\leq1\}$ of the $\TV$ norm is  relatively cheap computationally -- it reduces to solving a specific maximum flow problem.
It should be mentioned here that the relation between flow problems and TV-based denoising (problem (\ref{prob21}) with $\cA=I$) is well known and is utilized in many algorithms, see \cite{Goldfarb:2009} and references therein. While we have no doubt that the simple fact stated Lemma \ref{lemflow} below is well-known, for reader convenience we present here in detail the reduction mechanism.
\par
Consider the network (the oriented graph) $G$ with $n^2$ nodes $[i;j]\in\Gamma_{n,n}$ and $2n(n-1)$ arcs as follows: the first $n(n-1)$ arcs are of the form $([i+1;j],[i;j])$, $0\leq i<n-1$, $0\leq j<n$, the next $n(n-1)$ arcs are $([i;j+1],[i;j])$, $0\leq i<n$, $0\leq j<n-1$, and the remaining $2n(n-1)$ arcs (let us call them {\sl backward} arcs) are the inverses of the just defined $2n(n-1)$ {\sl forward} arcs. Let $\cE$ be the set of arcs of our network, and let us equip all the arcs with unit capacities. Let us treat vectors from $E=M^n_0$ as vectors of external supplies for our network; note that the entries of these vectors sum to zero, as required from external supply. Now, given a nonzero vector $\eta\in M^n_0$, let us consider the network flow problem where we seek for the largest multiple $s\eta$ of $\eta$ which, considered as the vector of external supplies in our network, results in a feasible capacitated network flow problem. The problem in question reads
\begin{equation}\label{network}
s_*=\max_{s,r}\left\{s: Pr=s\eta, \;0\leq r\leq {\mathbf{e}}\right\},
\end{equation}
where $P$ is the incidence matrix of our network\footnote{that is, the rows of $P$ are indexed by the nodes, the columns are indexed by the arcs, and in the column indexed by an arc $\gamma$ there are exactly two nonzero entries: entry 1 in the row indexed by the starting node of $\gamma$, and entry $-1$ in the row indexed by the terminal node of $\gamma$.} and $\mathbf{e}$ is the all-ones vector. Now, problem (\ref{network}) clearly is feasible, and its feasible set is bounded due to $\eta\neq0$, so that the problem is solvable.  Due to its network structure, this LP program can be solved reasonably fast even in the large scale case (say, when $n=512$ or $n=1024$, which already is of interest for actual imaging). Further, an intelligent network flow solver as applied to (\ref{network}) will return not only the optimal $s=s_*$ and the corresponding flow, but also the dual information, in particular, the optimal vector $z$ of Lagrange multipliers for the linear equality constraints $Pr-s\eta=0$. Let $\bar{z}$ be obtained by subtracting from the entries of $z$ their mean; since the entries of $z$ are indexed by the nodes, $\bar{z}$ can be naturally interpreted as a zero mean image. It turns out that this image is nonzero, and the vector $x[\eta]=-\bar{z}/\TV(\bar{z})$ is nothing than a desired minimizer of $\langle\eta,\cdot\rangle$ on $\CTV$:
\begin{lemma}\label{lemflow} Let $\eta$ be a nonzero image with zero mean. Then (\ref{network}) is solvable with positive optimal value, and the image $x[\eta]$, as defined above, is well defined and is a maximizer of $\langle \eta,\cdot\rangle$ on $\CTV$.
\end{lemma}
\paragraph{Bounding $L_f$.} When applying {\CG} algorithms to the TV-based problems (\ref{prob1}), (\ref{prob21}) with $E=M^n_0$ and $f(x)=\phi(\cA x - b)$, the efficiency estimates depend linearly on the associated quantity $L_f$, which, in turn, is readily given by the norm $\|\cA\|_{\tv(\cdot),\pi(\cdot)}$ of the mapping $x\mapsto \cA x$, see the end of Section \ref{sec:prb}. Observe that in typical applications $\cA$ is a simple operator (e.g., the discrete convolution), so that when restricting ourselves to the case when $\pi(\cdot)$ is $\|\cdot\|_2$  (quadratic fit),  it is easy to find a tight upper bound on $\|\cA\|_{\|\cdot\|_2,\|\cdot\|_2}$. To convert this bound into an upper bound on $\|\cA\|_{\tv(\cdot),\|\cdot\|_2}$, we need to estimate the quantity
$$
Q_n=\max_x\{\|x\|_2:x\in M^n_0,\TV(x)\leq1\}.
$$
Bounding $Q_n$ is not a completely trivial question, and the answer is as follows:
\begin{proposition}\label{laplace}
$Q_n$ is nearly constant, specifically, $Q_n\leq O(1)\sqrt{\ln(n)}$ with a properly selected absolute constant $O(1)$.
\end{proposition}
Note that the result of Proposition \ref{laplace} is  in sharp contrast with one-dimensional case, where the natural analogy of $Q_n$ grows with $n$ as $\sqrt{n}$. We do not know whether it is possible to replace in Proposition \ref{laplace} $O(1)\sqrt{\ln(n)}$ with $O(1)$, as suggested by Sobolev's inequalities\footnote{From the Sobolev embedding theorem it follows that for a smooth function $f(x,y)$ on the unit square one has $\|f\|_{L_2}\leq O(1)\|\nabla f\|_1,$ $\|\nabla f\|_1:=\|f'_x\|_1+\|f'_y\|_1$, provided that $f$ has zero mean. Denoting by $f^n$ the restriction of the function onto a $n\times n$ regular grid in the square, we conclude that $\|f^n\|_2/\TV(f^n)\to \|f\|_{L_2}/\|\nabla f\|_1\leq O(1)$ as $n\to\infty$. Note that the convergence in question takes place only in the 2-dimensional case.}. Note that on inspection of the proof, Proposition extends to the case of $d$-dimensional, $d>2$, images with zero mean, in which case $Q_n\leq C(d)$ with appropriately chosen $C(d)$.

\section{Numerical examples}\label{sec:numex}
We present here some very preliminary simulation results.
\subsection{{\CG} for parametric optimization: sparse matrix completion problem}
The
goal of the first series of our experiments is to illustrate how the performance and requirements of {\CG} algorithm for parametric optimization, when applied to the matrix completion problem \cite{care2009exact}, scale with  problem size. Specifically, we apply the algorithm of Section \ref{subsec:posec} to the problem of nuclear norm minimization
\be
\min\|\sigma(x)\|_1,\;\;\mbox{subject to}\;\sum_{(i,j)\in \Omega}(y_{ij}-x_{ij})^2\le \delta,
\ee{exactmc}
where $\sigma(x)$ is the singular spectrum of a $p\times q$ matrix $x$. In our experiments, the set $\Omega$ of observed entries $(i,j)\in \{1,...,p\}\times\{1,...,q\}$ of cardinality $m\ll pq$ was selected at random.
\par
 Note that the the implementation of the {\CG}M is especially simple for the problem \rf{exactmc} -- at each method's iteration it requires solving a simple quadratic problem with dimension of the decision variable which does not exceed the iteration count. This allows to implement efficiently the ``full memory'' version of {\CG}M ({\CG} algorithms with memory) \rf{bundlex}, \rf{bundleaux}, in which the set $X_t$ contains  $x_t$ and all the points $x^+_\tau$ for $1\le \tau\le t$.
\par
We compare the performance of  {\CG}M algorithms and of a ``memoryless'' version of the {\CG}. To this end we have conducted the following experiment:
\begin{enumerate}
\item For matrix sizes $p,q\in [1,2,4,8,16,32]\times 10^3$  we generate $n=10$ sparse $p\times q$ matrices $y$ with density $d=0.1$ of non-vanishing entries as follows: we generate $p\times r$ matrix $U$ and $q\times r$ matrix $V$ with independent Gaussian entries $u_{ij}\sim \cN(0,m^{-1}),\; v_{ij}\sim \cN(0,n^{-1})$, and a $r\times r$ diagonal matrix $D={\rm diag}[d_1,...,d_r]$ with $d_i$ drawn independently from a uniform distribution on $[0,1]$. The non-vanishing entries of the sparse observation matrix $y$  are obtained by sampling at random with probability $d$ the entries of $x^*=UDV^T$, so that for every $i,j$, $y_{ij}$ is, independently over $i,j$, set to $x^*_{ij}$ with probability $d$ and to $0$ with probability $1-d$. Thus, the number of non-vanishing entries of  $y$ is approximately $m=dpq$. This procedure results in $m\sim 10^5$ for the smallest matrices $y$ ($1000\times 1000$), and in $m\sim 10^8$ for  the largest matrices ($32000\times 32000$).
\item We apply to parametric optimization problem \rf{exactmc} MATLAB implementations of the {\CG}M with memory parameter $M=1$ (``memoryless'' {\CG}), {\CG}M with $M=5$ and full memory {\CG}M.
 The parameter $\delta$ of \rf{exactmc} is chosen to be $\delta=0.001 \|y\|^2_{\rm f}$ (here $\|y\|_{\rm f}=\left(\sum_{i,j} y^2_{ij}\right)^{1/2}$ stands for the Frobenius norm of $y$). The optimization algorithm is tuned to the {\em relative accuracy} $\varepsilon=1/4$, what means that it outputs an $\epsilon$-solution  $\widehat{x}$ to \rf{exactmc}, in the sense of \rf{goalgbspp}, with absolute accuracy $\epsilon=\delta\varepsilon$.
    \end{enumerate}
 For each algorithm (memoryless {\CG}, {\CG}M with memory $M=5$ and full memory {\CG}M) we present in table \ref{tab:table1} the average, over algorithm's runs on the (common for all algorithms) sample of $n=10$ matrices $y$ we have generated, 1) total number of iterations $N_{\rm it}$ necessary to produce an $\epsilon$-solution (it upper-bounds the rank of the resulting $\epsilon$-solutuion), 2) CPU time in seconds $T_{\rm cpu}$ and 3) MATLAB memory usage in megabytes $S_{\rm mem}$. This experiment was conducted on a {Dell Latitude 6430 laptop} equipped with  Intel Core i7-3720QM CPU@2.60GHz and 16GB of RAM. Because of high memory requirements in our implementation of the full memory {\CG}M, this method was unable to complete the computation for the two largest matrix sizes.
 \par
 We can make the following observation regarding the results summarized in table \ref{tab:table1}: {\CG} algorithm with memory consistently outperforms the standard -- memoryless -- version of {\CG}. The full memory {\CG}M requires the smallest number of iteration to  produce an $\epsilon$-solution, which is of the smallest rank, as a result. On the other hand, the memory requirements of the full memory {\CG}M become prohibitive (at least, for the computer we used for this experiment and MATLAB implementation of the memory heap) for large matrices. On the other hand, a {\CG}M with memory $M=5$ appears to be a reasonable compromise in terms of numerical efficiency and memory demand.
\begin{table}
\begin{center}
{\footnotesize\begin{tabular}{|c||c|c|c||c|c|c||c|c|c|}
\hline
Matrix size& \multicolumn{3}{|c||}{Memory-less {\CG}}&\multicolumn{3}{c||}{{\CG}M with memory $M=5$}&\multicolumn{3}{c|}{Full memory {\CG}}\\
\cline{2-4}\cline{5-7}\cline{8-10}
$p\times q$&$N_{\rm it}$&$T_{\rm cpu}$&$S_{\rm mem}$&$N_{\rm it}$&$T_{\rm cpu}$&$S_{\rm mem}$&$N_{\rm it}$&$T_{\rm cpu}$&$S_{\rm mem}$\\
\hline\hline
$1000\times1000$&271.6&9.35&17.11&149.7&5.01&17.63&78.4&4.71&78.98\\
\hline
$1000\times2000$& 292.1 &  12.14&   31.67&  162.8&    7.76&   32.57&   93.5&   10.89&  156.22\\
\hline
$2000\times2000$& 246.8&   17.01&   54.45&  139.1&   11.19&   61.57&   71.9&   13.31&  248.13\\\hline
$2000\times4000$& 259.3 &  33.94&  105.09&  152.3&   24.50&  120.22&   57.7&   25.54&  410.02\\\hline
$4000\times4000$& 321.8 &   79.20&   207.26&   162.9&    50.59&   235.59&    74.6&    93.22&  1014.7\\\hline
$4000\times8000$&360.1&   169.8&   399.16&   147.3&    88.81   &464.68&    63.3&   135.6  &1766.4\\\hline
$8000\times8000$&323.4 &     302.8&      754.46&   111.8&      134.1&      905.98&  53.6&      191.3&        3061.5\\\hline
$8000\times16000$&324.1 &     614.3&      1485.6&    118.2&      286.5&      1800.7&    50.5  &    329.4&      5826.7\\\hline
$16000\times16000$&258.7 &     995.4&      2898.5&    99.7&      495.5&      3577.8&    50.8  &    595.2 &     11696\\\hline
$16000\times32000$&276.7&      2572&      5721.7&      70.3&      859.2&      7109.0 &NA    & NA& NA\\\hline
$32000\times32000$&305.4&      5028&      11352&      57.6&      2541&      14186&NA &NA & NA\\
\hline\end{tabular}}
\caption{\label{tab:table1} memoryless {\CG} vs. {\CG}M with memory $M=5$ vs. full memory {\CG}M. $N_{\rm it}$: total number of method iterations; $T_{\rm cpu}$: CPU usage (sec), and $S_{\rm mem}$: memory usage (MB) reported by MATLAB.}
\end{center}
\end{table}


\subsection{{\CG} for composite optimization: multi-class classification with nuclear-norm regularization}
We present here an empirical study of the {\CG} algorithm for composite optimization
as applied to the machine learning problem of multi-class classification with nuclear-norm penalty. A brief description of the multi-class classification
 problem is as follows: we observe $N$ ``feature vectors'' $\xi_i\in\bR^q$, each belonging to exactly one of $p$ classes $C_1,...,C_p$. Each $\xi_i$ is augmented by its label $y_i\in\{1,...,p\}$ indicating to which class $\xi_i$ belongs. Our goal is to build a classifier capable to predict the class to which a new feature vector $\xi$ belongs. This classifier is given by a $p\times q$ matrix $x$ according to the following rule: given $\xi$, we compute the $p$-dimensional vector $x\xi$ and take, as the guessed class of $\xi$, the index of the largest entry in this vector.\par
In some cases (see \cite{dhm:2012,hadopaduma:2012}), when, for instance, one is dealing with a large number of classes, there are good reasons ``to train the classifier'' --- to specify $x$ given the training sample $(\xi_i,y_i)$, $1\leq i\leq N$ --- as the optimal solution to the nuclear norm penalized minimization problem
\begin{equation}
\label{eq:nuclearmulticlassobj}
\Opt(\kappa)=\min_{x\in\bR^{p\times q}} \quad F_\kappa(x):=\overbrace{\frac{1}{N} \sum_{i=1}^{N} \log \left\{
\sum_{\ell=1}^q \exp \left((x_{\ell}^{T} -x_{y_i}^{T})\xi_i \right)
\right\}}^{f(x)}
+ \kappa \|\sigma(x)\|_1,
\end{equation}
where $x_\ell^T$ is the $\ell$-th row in $x$. \par
Below, we report on some experiments with this problem. Our goal was to compare
two versions of {\CG} for composite minimization: the memoryless version defined in (\ref{cgco}) and the version with memory
defined in (\ref{bundle3}). To solve the corresponding sub-problems, we used
the Center of Gravity method in the case of (\ref{cgco}) and the Ellipsoid method in the case of (\ref{bundle3})~\cite{Nesterov:2004,Nemirovski:Yudin:1983}.
In the version with memory we set $M=5$, as it appeared to be the best option from empirical evidence.
We have considered the following datasets:
\begin{enumerate}

\item \textbf{Simulated data}: for matrix of sizes $p,q\in 10^3 \times \{2^s\}_{s=1}^4$, we generate random matrices $x_\star = U S V$, with $p\times p$ factor $U$, $q\times q$ factor $V$, and diagonal $p\times q$ factor $S$  with random entries
sampled, independently of each other, from $\cN(0,p^{-1})$ (for $U$), $\cN(0,q^{-1})$ (for $V$), and the uniform distribution on $[0,1]$ (for diagonal entries in $S$).
We use $N=20 q$, with the feature vectors  $\xi_1,...,\xi_N$ sampled, independently of each other, from the distribution $\cN(0,I_q)$, and their labels $y_i$ being the indexes of the largest entries in the vectors
$x_\star\xi_i+\epsilon_i$, where $\epsilon_i\in\bR^p$ were sampled, independently of each other and of $\xi_1,...,\xi_N$, from $\cN(0,{1\over 2}I_p)$. The regularization parameter
$\kappa$ is set to $10^{-3}\Tr(x_\star x_\star^T)$.

\item \textbf{Real-world data}: we follow a setting similar to~\cite{hadopaduma:2012}.
We consider the Pascal ILSVRC2010 ImageNet dataset and focus on the ``Vertebrate-craniate'' subset, yielding  $1043$ classes, {with} $20$ examples per class.
The goal here is to train a multi-class classifier in order to be able to predict the class of each image (example) of the dataset.
Each example is converted to a $65536$-dimensional feature vector of unit $\ell_1$-norm using state-of-the-art visual descriptors known as Fisher
vector representation~\cite{hadopaduma:2012}. To summarize, we have $p=1043$, $q=65536$,  $N=20860$.
We set the regularization parameter to $\kappa = 10^{-4}$, which was found to result in the best predictive performance as estimated by cross-validation,
a standard procedure to set the hyper parameters in machine learning~\cite{Hastie:Tibshirani:Friedman:2008}.

\end{enumerate}

In both sets of experiments, the computations are terminated when the ``$\epsilon$-optimality conditions''
\begin{equation}
\label{co-cg-stop}
\begin{array}{rcl}
\|\sigma(f'(x_t))\|_\infty &\leq&\kappa + \epsilon \\
\langle f'(x_t), x_t\rangle + \kappa \|\sigma(x_t)\|_1 &\leq& \epsilon \|\sigma(x_t)\|_1 \\
\end{array}
\end{equation}
were met,
where $\| \sigma(\cdot) \|_{\infty}$ denotes the usual operator norm (the largest singular value). These conditions admit transparent interpretation as follows. For every $\bar{x}$, the function
\begin{equation*}
\phi_\kappa(x)=f(\bar{x})+\langle f'(\bar{x}),x-\bar{x}\rangle +\kappa \|\sigma(x)\|_1
\end{equation*}
underestimates $F_\kappa(x)$, see (\ref{eq:nuclearmulticlassobj}),
whence $\Opt(\kappa') \geq f(\bar{x})-\langle f'(\bar{x}),\bar{x}\rangle$ whenever
$\kappa'\geq \|\sigma(f'(\bar{x}))\|_\infty$. Thus, whenever $\bar{x}=x_t$ satisfies the first relation in (\ref{co-cg-stop}), we have $\Opt(\kappa+\epsilon)\geq f(x_t)-\langle f'(x_t),x_t\rangle$, whence \[F_\kappa(x_t)-\Opt(\kappa+\epsilon)\leq \langle f'(x_t), x_t\rangle + \kappa \|\sigma(x_t)\|_1.\] We see that (\ref{co-cg-stop}) ensures that
$F_\kappa(x_t)-\Opt(\kappa+\epsilon)\leq \epsilon \; \|\sigma(x_t)\|_1 $,
which, for small $\epsilon$, is a reasonable substitute for the actually desired termination when $F_\kappa(x_t)-\Opt(\kappa)$ becomes small.  In our experiments, we use $\epsilon=0.001$.

In table \ref{tab:table2} for each algorithm (memoryless {\CG}, {\CG}M with memory $M=5$) we present the average, over 20 collections of simulated data coming from 20 realizations of $x_\star$, of: 1) total number of iterations $N_{\rm it}$ necessary to produce an $\epsilon$-solution, 2) CPU time in seconds $T_{\rm cpu}$. The last row of the table corresponds to the real-world data. Experiments were conducted on a Dell R905 server equipped with four six-core AMD Opteron 2.80GHz CPUs and 64GB of RAM. A maximum of 32GB of RAM was used for the computations.

We draw the following conclusions from table \ref{tab:table1}: {\CG} algorithm with memory routinely outperforms the standard -- memoryless -- version of {\CG}. However, there is a trade-off between the algorithm progress at each iteration
and the computational load of each iteration. Note that, for large $M$, solving the sub-problem (\ref{bundle3}) can be challenging.

\begin{table}
\begin{center}
{\footnotesize\begin{tabular}{|c||c|c|c||c|c|c|c||}
\hline
Matrix size& \multicolumn{3}{|c||}{Memory-less {\CG}}&\multicolumn{3}{c||}{{\CG}M with memory $M=5$}\\
\cline{2-4}\cline{5-7}
$p\times q$&$N_{\rm it}$&$T_{\rm cpu}$&$S_{\rm mem}$&$N_{\rm it}$&$T_{\rm cpu}$&$S_{\rm mem}$\\
\hline\hline
$2000\times2000$&172.9&$349.7$&134.4&99.70&125.1 &174.1\\
$4000\times4000$&153.4&$ 1035 $&541.8&88.2&$575.2$ &704.1\\
$8000\times8000$&195.3&$ 2755$&2169&120.4&$1284$ &2819\\
$16000\times16000$&230.2&$ 6585$&8901&134.3&$3413$ &11550\\
$32000\times32000$&271.4&$26370$&30300&140.4&$17340$ & 30500\\ \hline
\hline
$1043\times65536$&183 &$2101$&2087&111&$925.34$ &2709\\
\hline\end{tabular}}
\caption{\label{tab:table2} memoryless {\CG} vs. {\CG}M with memory $M=5$. $N_{\rm it}$: total number of method iterations; $T_{\rm cpu}$: CPU usage (sec) reported by MATLAB.
}
\end{center}
\end{table}


\subsection{{\CG} for composite optimization: TV-regularized image reconstruction}
Here we report on experiments with CO{\CG}M as applied to TV-regularized image reconstruction. Our problem of interest is of the form (\ref{prob21}) with quadratic $f$, namely, the problem
\begin{equation}\label{prob21TV}
\min_{x\in M_0^n} \phi_\kappa(x):=\underbrace{\frac{1}{2}\|P\cA x - Pb\|_2^2}_{f(x)} +\kappa \TV(x);
\end{equation}
for notation, see section \ref{sect:TV}.
\paragraph{Test problems.} In our experiments, the mapping $x\mapsto \cA x$ is defined as follows: we zero-pad $x$ to extend it from $\Gamma_{n,n}$ to get a finitely supported function on ${\mathbf{Z}}^2$, then convolve this function with a finitely supported kernel $\alpha(\cdot)$, and restrict the result onto $\Gamma_{n,n}$. The observations $b\in M^n$ were generated at random according to
\begin{equation}\label{observTV}
b_{ij}=(\cA x)_{ij}+\sigma\|x\|_\infty\xi_{ij},\,\,\xi_{ij}\sim\cN(0,1),\,\,1\leq i,j\leq n,
\end{equation}
with mutually independent $\xi_{ij}$. The relative noise intensity $\sigma>0$, same as the convolution kernel $\alpha(\cdot)$, are parameters of the setup of an experiment.
\paragraph{The algorithm.} We used  the CO{\CG} with memory, described in section \ref{sec:cocg}; we implemented the options listed in {\bf A} -- {\bf C} at the end of the section. Specifically,
\begin{enumerate}
\item We use the updating rule (\ref{bundledoubleprime}) with $Z_t$ evolving in time exactly as explained in item {\bf C}: the set $Z_t$ is obtained from $Z_{t-1}$ by adding the points $z_t=[x_t;\TV(x_t)]$, $\widehat{z}_t=[x[\nabla f(x_t)];1]$ and $z_t^\prime=[\nabla f(x_t);\TV(\nabla f(x_t))]$, and deleting from the resulting set, if necessary, some ``old'' points, selected according to the rule ``first in -- first out,'' to keep the cardinality of $Z_t$ not to exceed a given $M\geq 3$ (in our experiments we use $M=48$). This scheme is initialized with $Z_0=\emptyset$, $z_1=[0;0]$.
\item We use every run of the algorithm to obtain a set of approximate solutions to (\ref{prob21TV}) associated with various values of the penalty parameter $\kappa$, as explained in {\bf B} at the end of section \ref{sec:cocg}. Precisely, when solving (\ref{prob21TV}) for a given value of $\kappa$ (in the sequel, we refer to it as to the {\sl working value}, denoted $\kappa_{\rm w}$), we also compute approximate solutions $x_\kappa(\kappa')$ to the problems with the values $\kappa'$ of the penalty, for $\kappa'=\kappa\gamma$, with $\gamma$ running through a given finite subset $G\ni 1$  of the positive ray. In our experiments, we used the 25-point grid $G=\{\gamma=2^{\ell/4}\}_{\ell=-12}^{12}$.
\end{enumerate}
\par
The {LO oracle} for the TV norm on $M_0^n$ utilized in CO{\CG}M was the one described in Lemma \ref{lemflow}; the associated flow problem (\ref{network}) was solved by the commercial interior point LP solver  {\tt mosekopt} version 6 \cite{Mosek}. Surprisingly, in our application this ``general purpose'' interior point LP solver was by orders of magnitude faster than all dedicated network flow algorithms we have tried, including simplex-type network versions of  {\tt mosekopt}
and {\tt CPLEX}. With our solver, it becomes possible to replace in (\ref{network}) every pair of opposite to each other arcs with a single arc, passing from the bounds $0\leq r\leq {\mathbf{e}}$ on the flows in the arcs to the bounds $-{\mathbf{e}}\leq r\leq {\mathbf{e}}$.
\par
The {\sl termination criterion} we use relies upon the fact that in CO{\CG}M the (nonnegative) objective decreases along the iterates: we terminate a run when the progress in terms of the objective becomes small, namely, when the condition
$$
\phi_\kappa(x_{t-1})-\phi_\kappa(x_t)\leq \epsilon\max[\phi_\kappa(x_{t-1}),\delta\phi_\kappa(0)]
$$
is satisfied. Here $\epsilon$ and $\delta$ are small tolerances (we used $\epsilon=0.005$ and $\delta=0.01$).
\paragraph{Organization of the experiments.} In each experiment we select  a ``true image'' $x^*\in M^n$,  a kernel $\alpha(\cdot)$ and a (relative) noise intensity $\sigma$. Then we generate a related observation $b$, thus ending up with a particular instance of (\ref{prob21TV}). This instance is solved by  the outlined algorithm for working values $\kappa_{\rm w}$ of $\kappa$ taken from the set $G^+=\{\gamma=2^{\ell/4}\}_{\ell=-\infty}^\infty$, with the initial working value, selected in pilot runs, of the penalty underestimating the best -- resulting in the best recovery -- penalty.
\par
\begin{figure}
$$
\begin{array}{cc}
\resizebox{210pt}{120pt}{\includegraphics{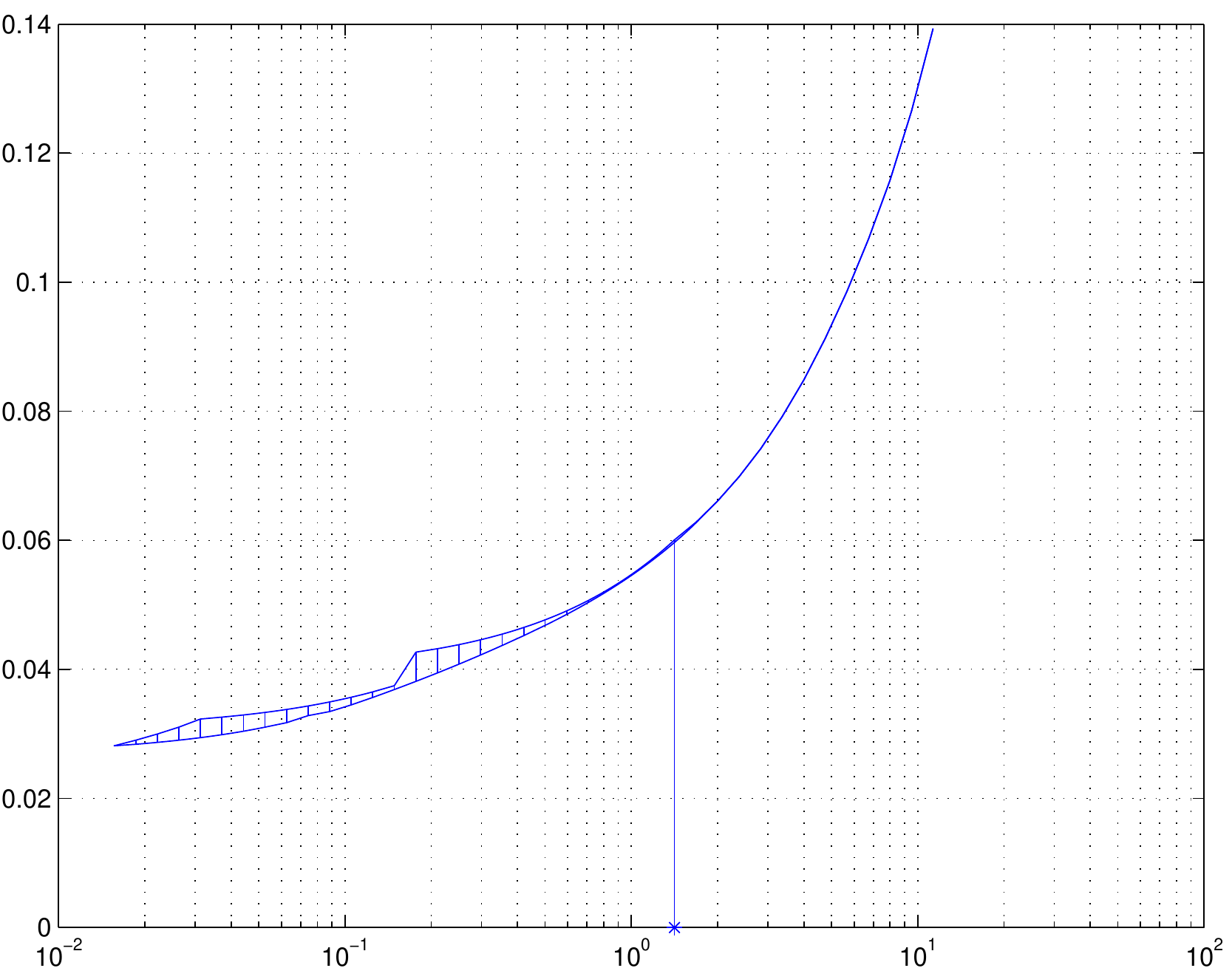}}&
\resizebox{210pt}{120pt}{\includegraphics{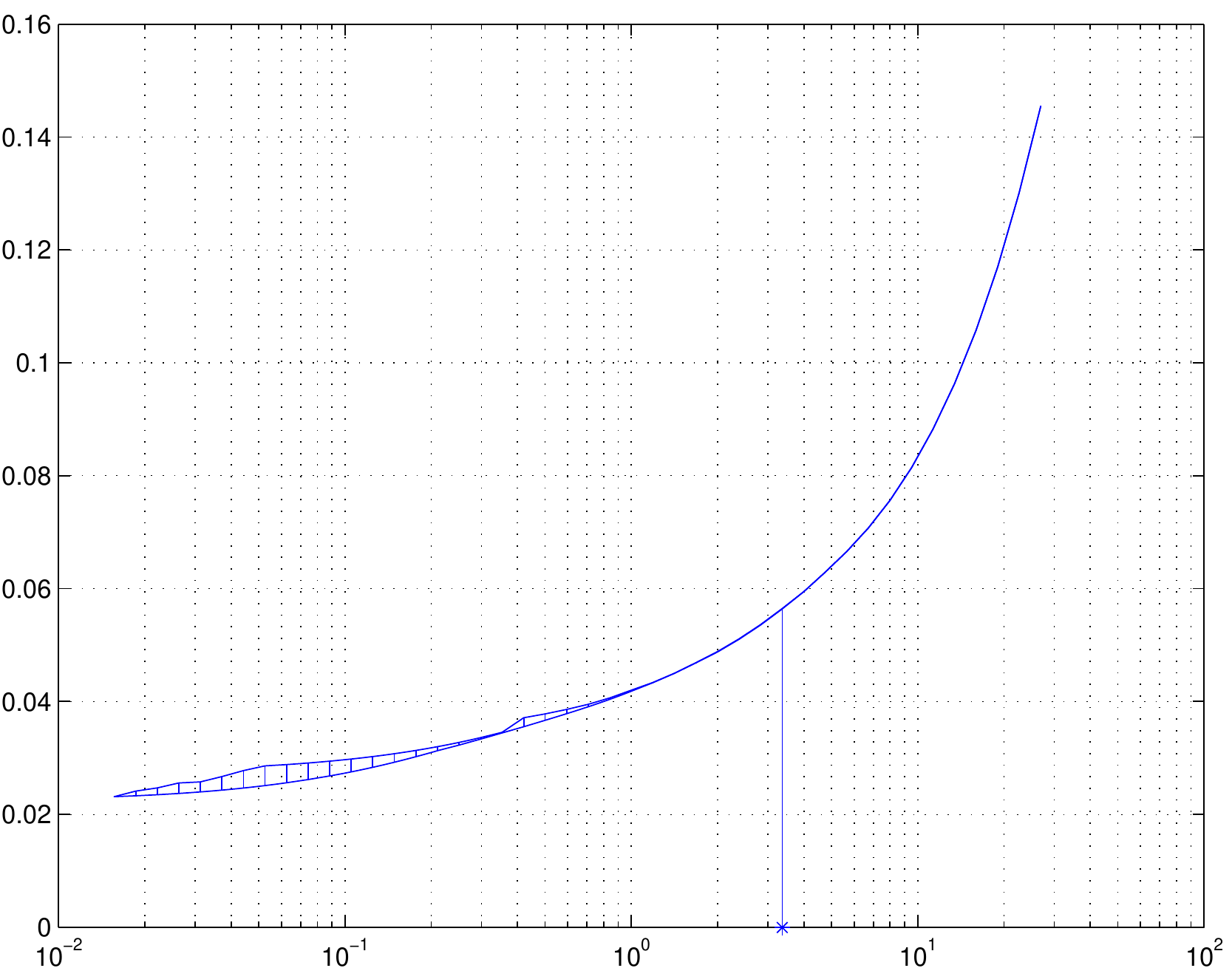}}\\
\hbox{(a)}&\hbox{(b)}\\
\resizebox{210pt}{120pt}{\includegraphics{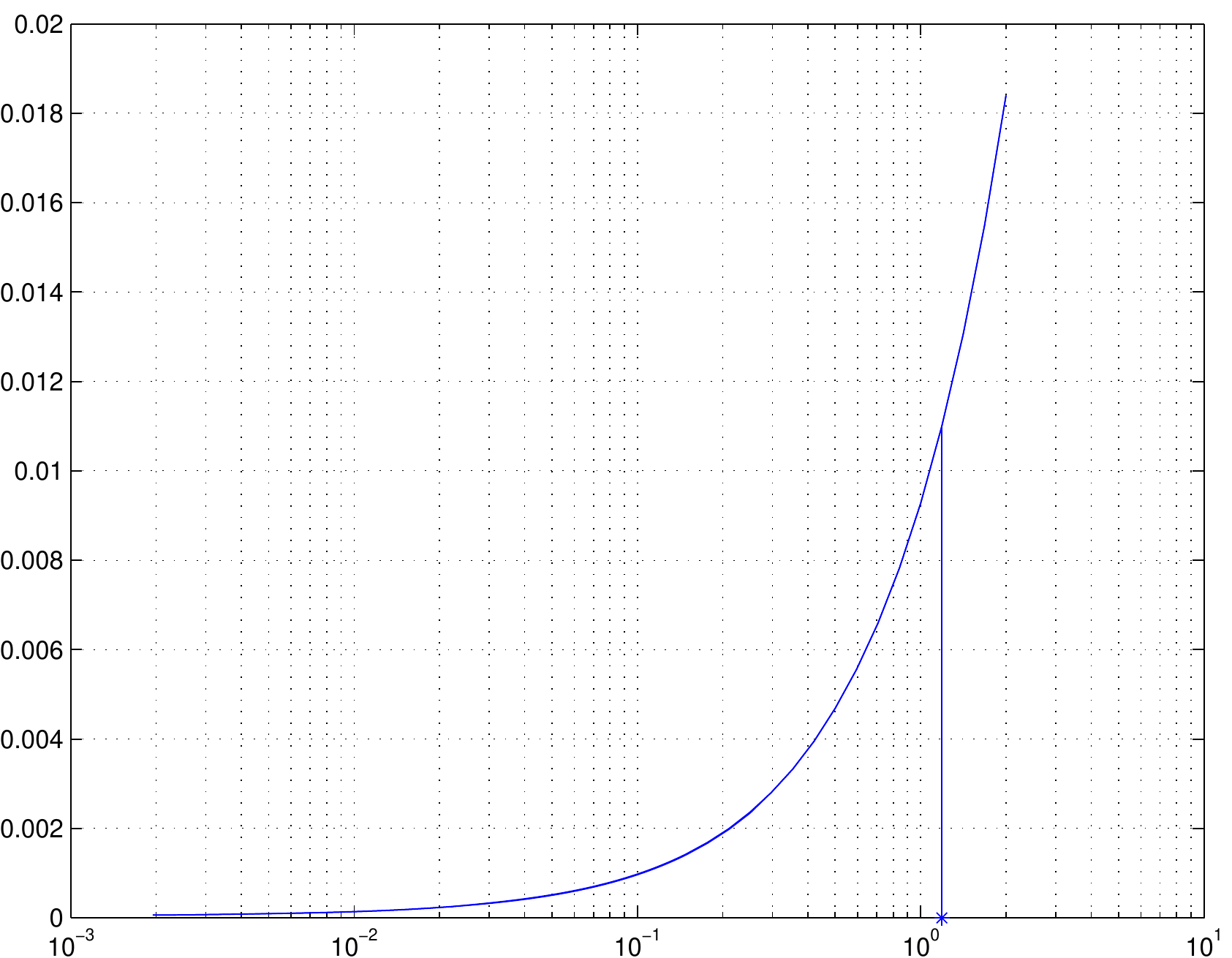}}&\resizebox{210pt}{120pt}{\includegraphics{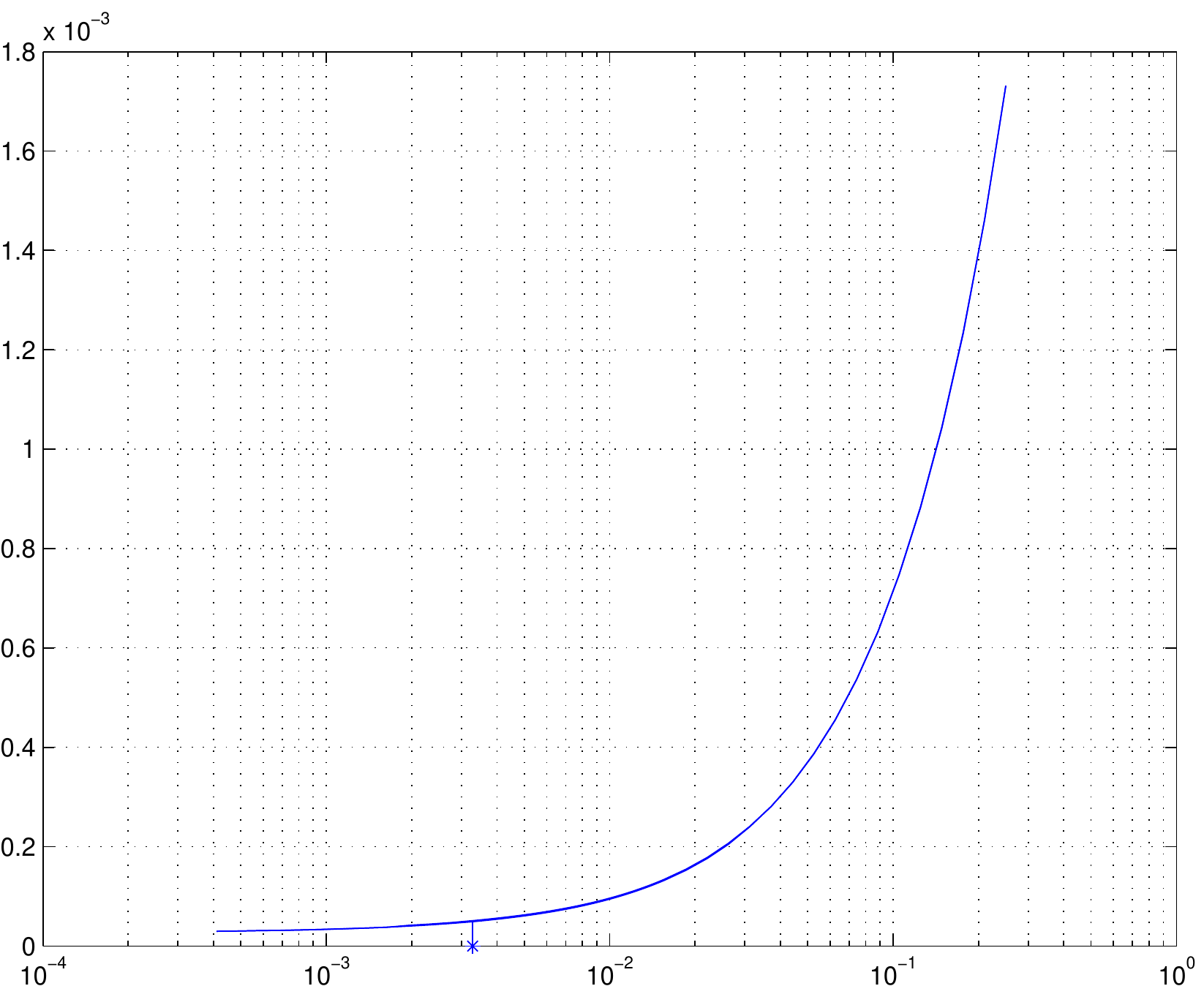}}\\
\hbox{(c)}&\hbox{(d)}\\
\multicolumn{2}{c}{\resizebox{210pt}{120pt}{\includegraphics{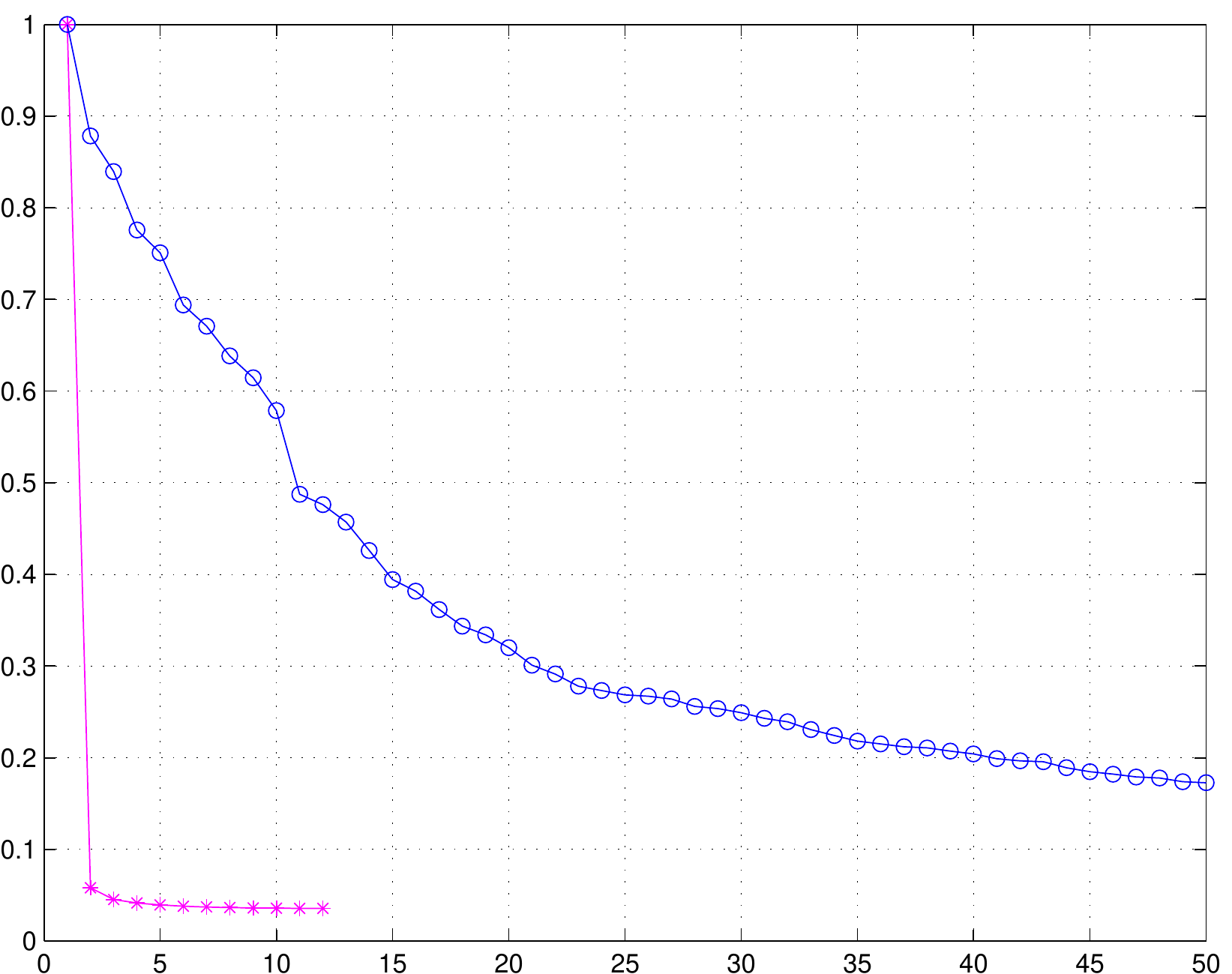}}}\\
\multicolumn{2}{c}{\hbox{(e)}}\\
\end{array}
$$
\caption{\small\label{figTV} (a) -- (d): lower and upper envelopes of $\{\phi_{\kappa}(x_{\kappa_{\rm w}}(\kappa)):\kappa_{\rm w}\in G\}$ vs. $\kappa$, experiments A -- D.
Asterisks on the $\kappa$-axes: penalties resulting in the smallest combined relative recovery errors.
(e): values of $\phi_{1/8}(x_t)$ vs. iteration number $t$ with $z_t^\prime$ included (asterisks) and not included (circles) into $Z_t$.}
\end{figure}

\begin{figure}
$$
\begin{array}{|ccc|}
\hline
\resizebox{150pt}{120pt}{\includegraphics{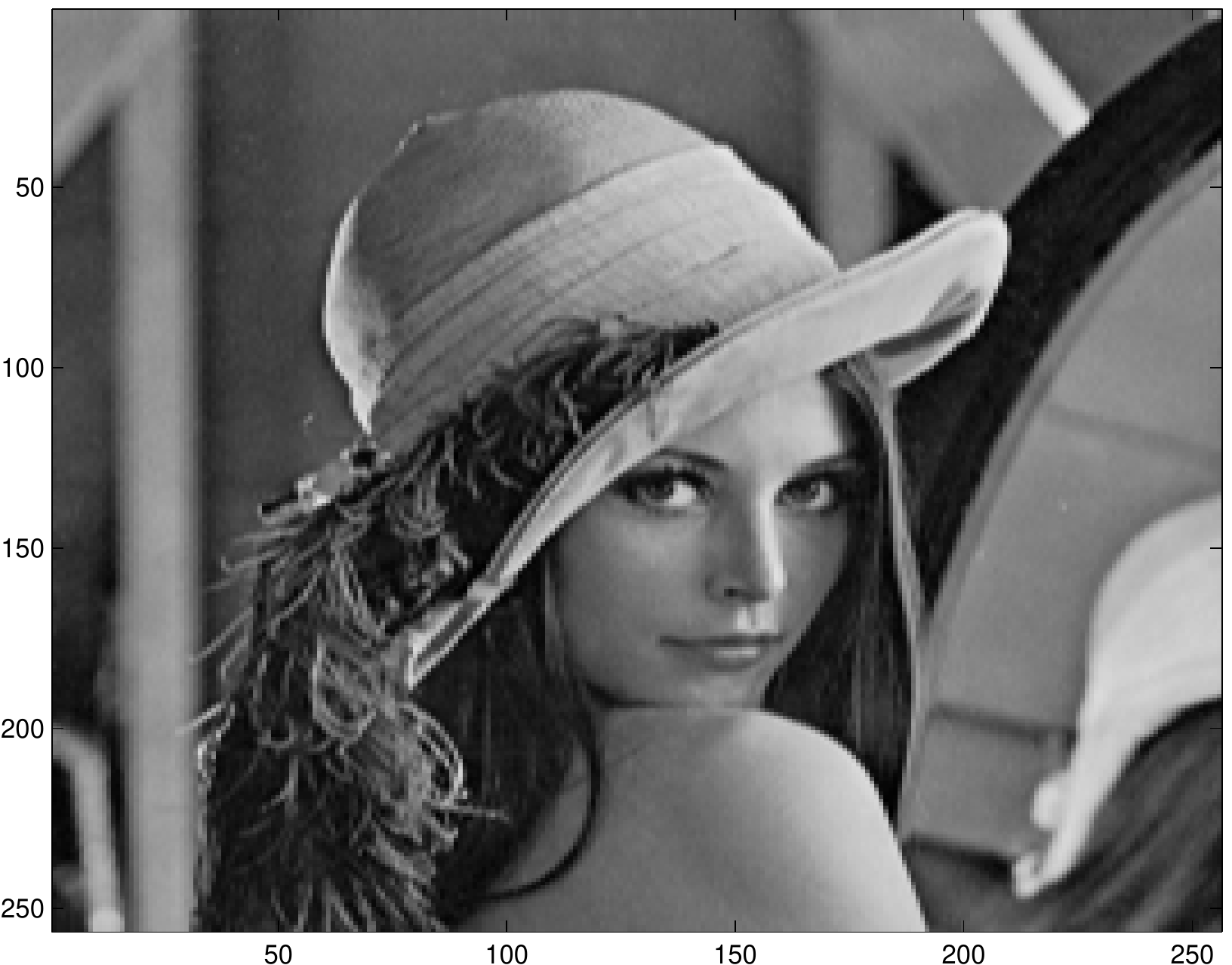}}&
\resizebox{150pt}{120pt}{\includegraphics{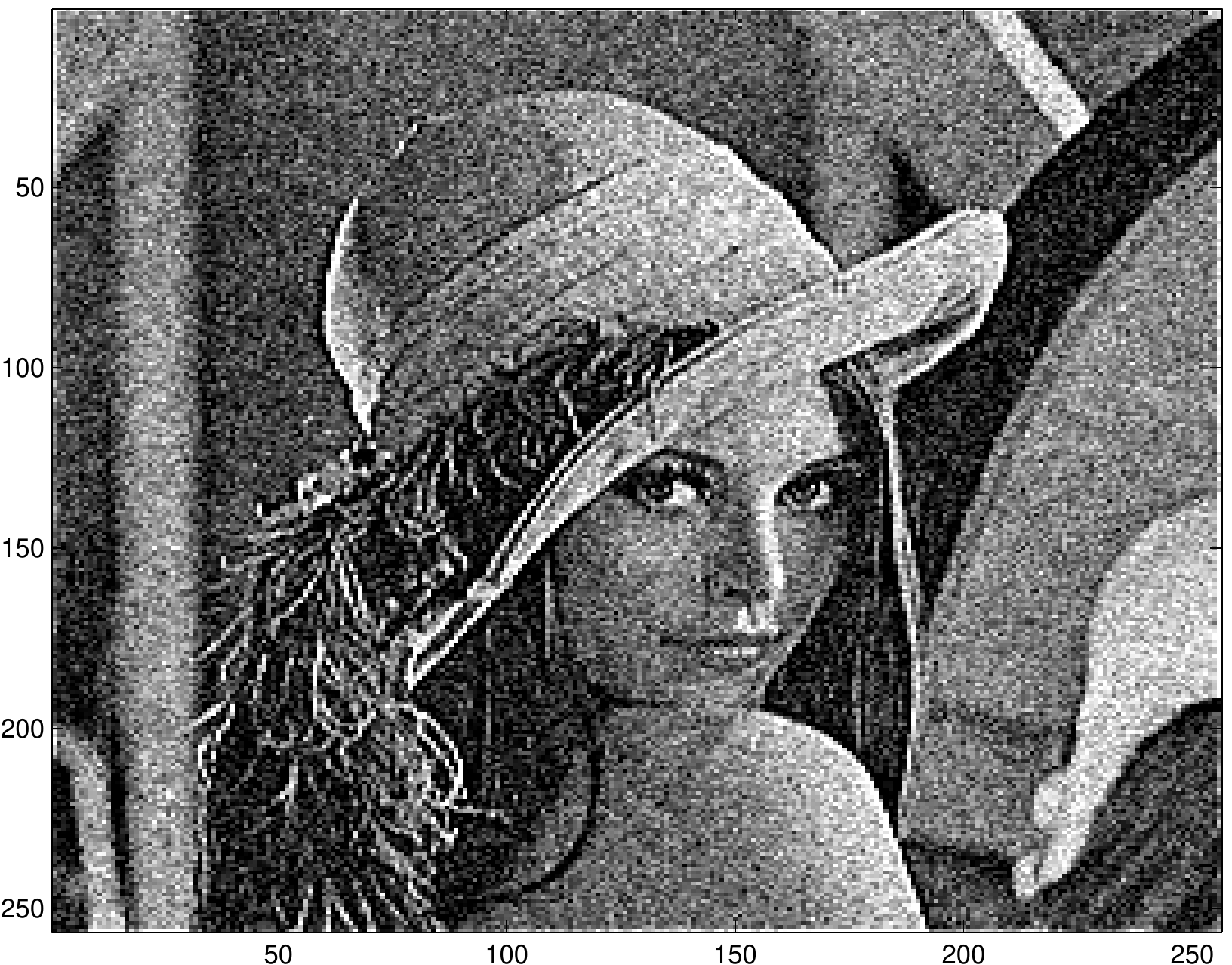}}&
\resizebox{150pt}{120pt}{\includegraphics{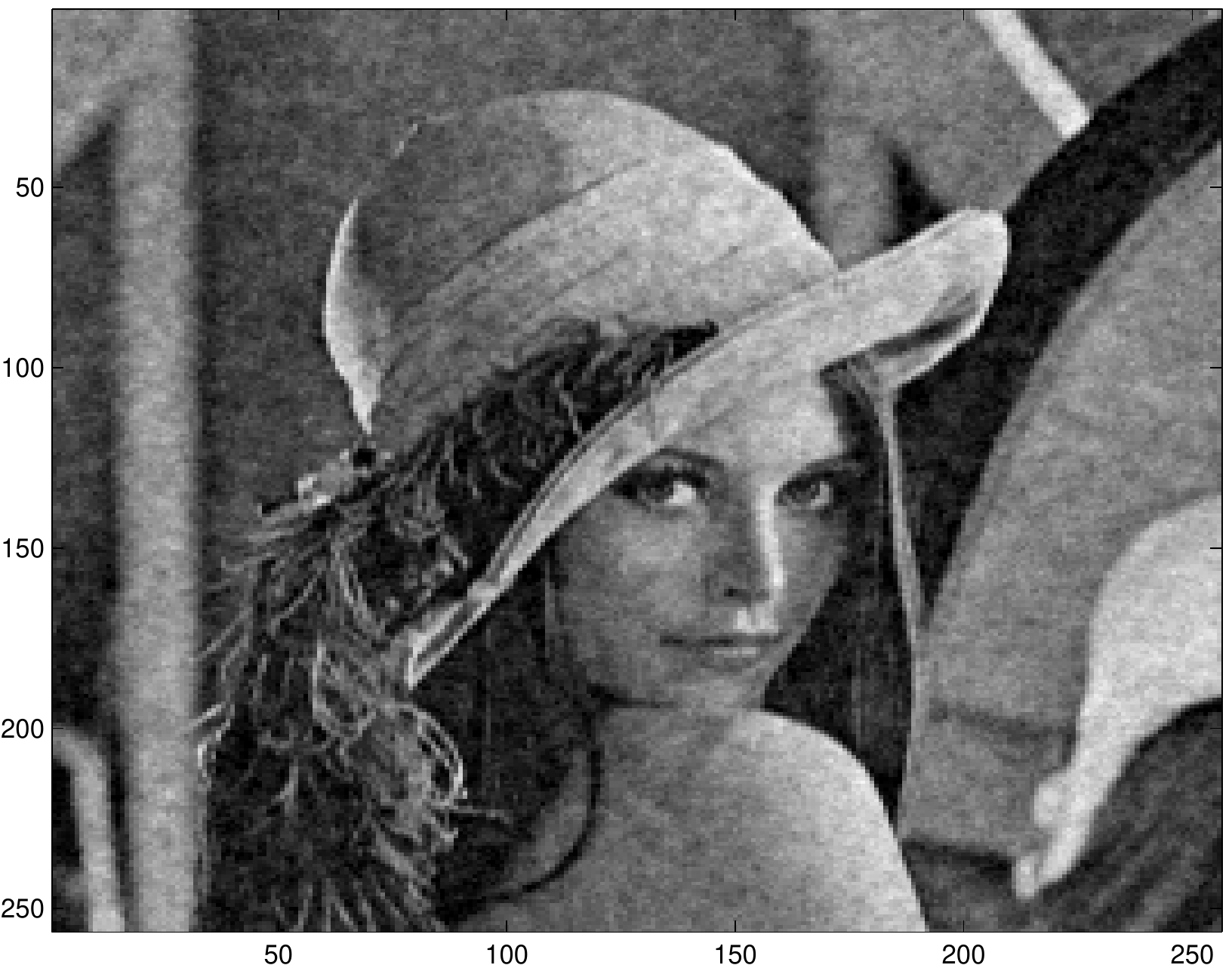}}\\
\hbox{C: True image}&\hbox{C: Observations, $\sigma=0.15$}&\hbox{C: TV recovery, $\kappa=0.250$}\\
\hline
\resizebox{150pt}{120pt}{\includegraphics{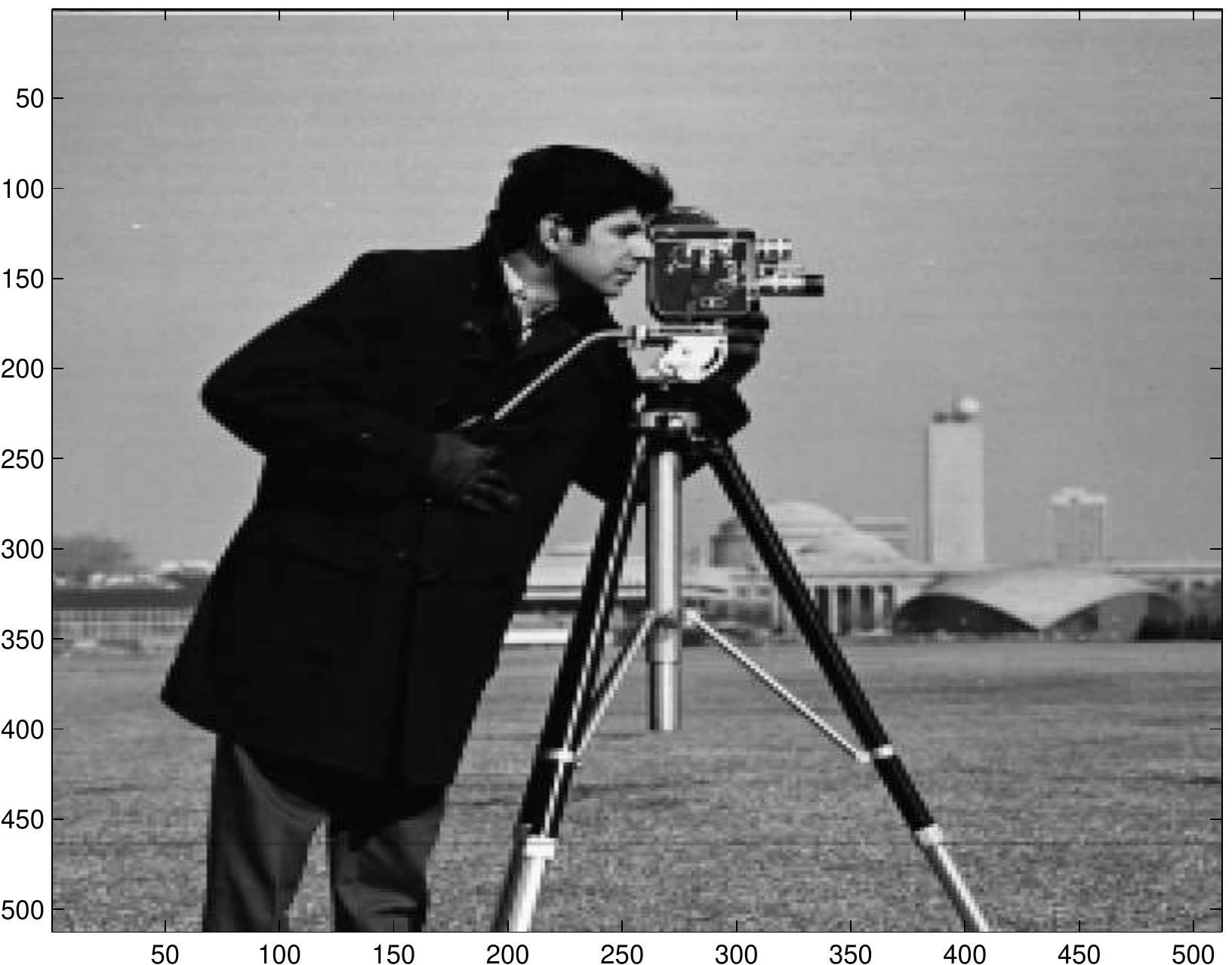}}&
\resizebox{150pt}{120pt}{\includegraphics{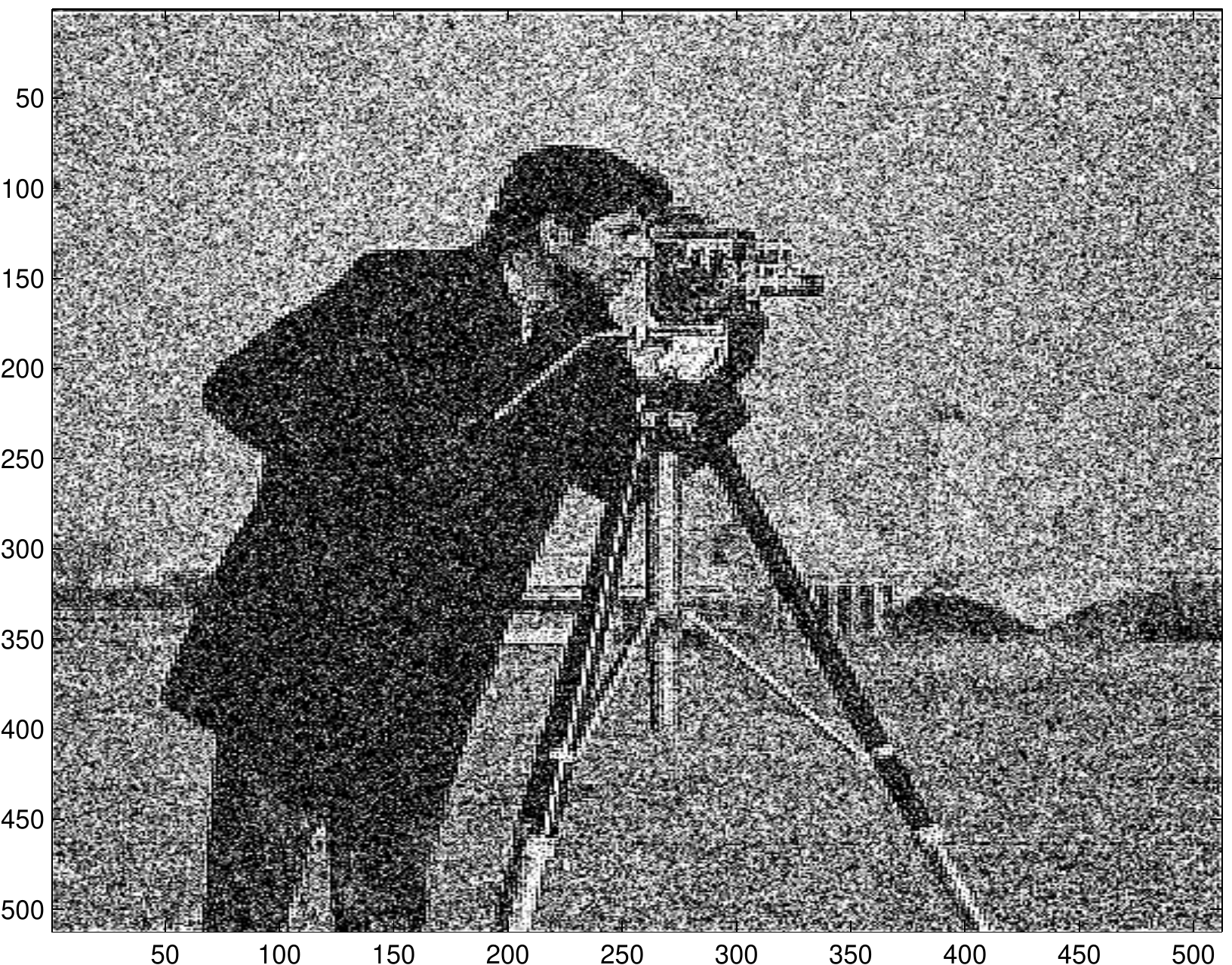}}&
\resizebox{150pt}{120pt}{\includegraphics{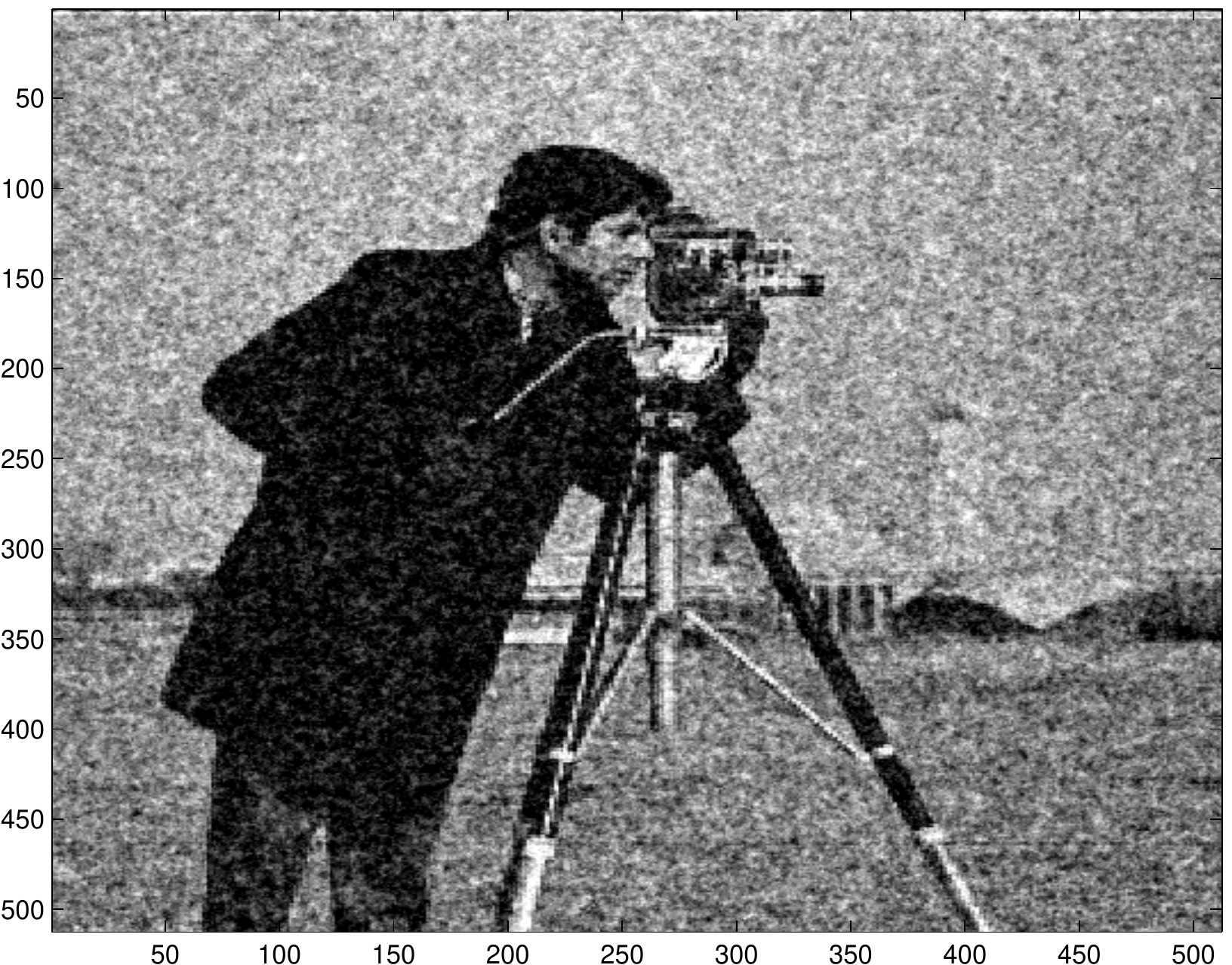}}\\
\hbox{D: True image}&\hbox{D: Observations, $\sigma=0.40$}&\hbox{D: TV recovery, $\kappa=0.00328$}\\
\hline
\end{array}
$$
\caption{\label{figTV2} Experiments C, D}
\end{figure}

 As explained above, a run of CO{\CG}M, the working value of the penalty being $\kappa_{\rm w}$, yields 25 approximate solutions to (\ref{prob21TV}) corresponding to $\kappa$ along the grid $\kappa_{\rm w}\cdot G$. These sets are fragments of the grid $G^+$, with the ratio of the consecutive grid points $2^{1/4}\approx 1.19$. For every approximate solution $x$  we compute its {\sl combined relative error} defined as
$$
\nu(x)=\left(\frac{\|\bar{x}-x^*\|_1\|\bar{x}-x^*\|_2\|\bar{x}-x^*\|_\infty}{\|x^*\|_1\|x^*\|_2\|x^*\|_\infty}\right)^{1/3};
$$
here $\bar{x}$ is the easily computable shift of $x$ by a constant image satisfying $\|\cA\bar{x}-b\|_2=\|P\cA x - Pb\|_2$.
 From run to run, we increase the working value of the penalty by the factor $2^{1/4}$, and terminate the experiment when  in four consecutive runs there was no progress in the combined relative error of the best solution found so far. Our primary goals are (a) {\sl to quantify the performance of the CO{\CG}M algorithm}, and (b) {\sl to understand by  which margin, in terms of $\phi_{\kappa}(\cdot)$, the ``byproduct'' approximate solutions yielded by the algorithm} (those which were obtained when solving (\ref{prob21TV}) with the working value of penalty different from $\kappa$) {\sl are worse than the ``direct'' approximate solution obtained for the working value $\kappa$ of the penalty.}
\paragraph{Test instances and results.} We present below the results of four experiments with two popular images; these results are fully consistent with those of other experiments we have conducted so far. The corresponding setups are presented in table \ref{tableTV1}. Table \ref{tableTV2} summarizes the performance data. Our comments are as follows.
\begin{itemize}
\item In accordance to the above observations, using ``large'' memory  (with the cardinality of $Z_t$ allowed to be as large as 48) and processing ``large'' number (25) of penalty values at every step are basically costless: at an iteration, the {\sl single} call to the {LO oracle} (which is a must for {\CG}) takes as much as $85\%$ of the iteration time.
\item The CO{\CG}M iteration count as presented in table \ref{tableTV2} is surprisingly low for an algorithm with sublinear $O(1/t)$ convergence, and the running time of the algorithm appears quite tolerable\footnote{For comparison: solving on the same platform problem (\ref{prob21TV}) corresponding to  Experiment A ($256\times256$ image) by the state-of-the-art commercial interior point solver {\tt mosekopt 6.0} took as much as 3,727 sec, and this -- for a {\sl single} value of the penalty (there is no clear way to get from a single run approximate solutions for a set of values of the penalty in this case).}
    \par
    Seemingly, the {\sl instrumental} factor here is that by reasons indicated in  {\bf C}, see the end of section \ref{sec:cocg}, we include into  $Z_t$ not only $z_t=[x_t;\TV(x_t)]$ and $\widehat{z}_t=[x[\nabla f(x_t)];1]$,  but also $z_t^\prime=[\nabla f(x_t);\TV(\nabla f(x_t))]$. To illustrate the difference, this is what happens in experiment A with the lowest (0.125) working value of penalty. With the outlined implementation, the run takes 12 iterations (111 sec), with the ratio
    $\phi_{1/8}(x_t)/\phi_{1/8}(x_1)$ reduced from 1 $(t=1)$ to 0.036 $(t=12)$. When $z_t^\prime$ is not included into $Z_t$, the termination criterion is not met even in 50 iterations (452 sec), the maximum iteration count we allow for a run, and in course of these 50 iterations the above ratio was reduced from 1 to  0.17, see plot e) on figure \ref{figTV}.
\item An attractive feature of the proposed approach is the possibility to extract from a single run, the working value of the penalty being $\kappa_{\rm w}$, suboptimal solutions $x_{\kappa_{\rm w}}(\kappa)$ for a bunch of instances of (\ref{prob21}) differing from each other by the values of the penalty $\kappa$. The related question is, of course, how good, in terms of the objective $\phi_{\kappa}(\cdot)$, are the ``byproduct'' suboptimal solutions $x_{\kappa_{\rm w}}(\kappa)$ as compared to those obtained when $\kappa$ is the working value of the penalty. In our experiments, the ``byproduct'' solutions were pretty good, as can be seen from plots (a) -- (c) on figure \ref{figTV}, where we see the upper and the lower envelopes of the values of $\phi_{\kappa}$ at the approximate solutions $x_{\kappa_{\rm w}}(\kappa)$ obtained from different working values $\kappa_{\rm w}$ of the penalty. In spite of the fact that in our experiments the ratios $\kappa/\kappa_{\rm w}$ could be as small as $1/8$ and as large as $8$, we see that these envelopes are pretty close to each other, and, as an additional bonus, are merely indistinguishable in a wide neighborhood of the best (resulting in the best recovery) value of the penalty (on the plots, this value is marked by asterisk).
\end{itemize}
Finally, we remark that in experiments A, B, where the mapping $\cA$ is heavily ill-conditioned (see table \ref{tableTV1}), TV regularization yields moderate (just about 25\%) improvement in the combined relative recovery error as compared to the one of the trivial recovery (``observations as they are''), in spite of the relatively low ($\sigma=0.05$) observation noise. In contrast to this, in the experiments C, D, where $\cA$ is well-conditioned, TV regularization reduces the error by $80\%$ in experiment C ($\sigma=0.15$) and by 72\% in experiment D ($\sigma=0.4$), see figure \ref{figTV2}.
\begin{table}
\begin{center}
{\footnotesize
\begin{tabular}{||c||c|c|c|c|c||}
\hline\hline
\# &Image&$n$&$\alpha(\cdot)$&$\hbox{\rm Cond}(\cA^*\cA)$&$\sigma$\\
\hline\hline
A&{\tt lenna}$^\dag$&256&{\tt fspecial('gaussian',7,1)$^\S$} (7$\times$7)&$\approx2.5e7$&0.05\\
\hline
B&{\tt cameraman}$^\ddag$&512&{\tt fspecial('gaussian',7,1)} (7$\times$7)&$\approx2.5e7$&0.05\\
\hline
C&{\tt lenna}&256&{\tt fspecial('unsharp')} (3$\times$3)&$\approx40$&0.15\\
\hline
D&{\tt cameraman}&512&{\tt fspecial('unsharp')} (3$\times$3)&$\approx40$&0.40\\
\hline\hline
\multicolumn{6}{l}{\footnotesize$^\dag$http://en.wikipedia.org/wiki/Lenna$\quad$$^\ddag$http://en.wikipedia.org/wiki/Camera\_operator}\\
\multicolumn{6}{l}{\footnotesize$^\S$http://www.mathworks.com/help/images/ref/fspecial.html}\\
\end{tabular}}
\caption{\label{tableTV1} Setups of the experiments. }
\end{center}
\end{table}

\begin{table}
\begin{center}
{\footnotesize
\begin{tabular}{||c||c||c||c|c|c||c|c||c||}
\cline{4-9}
\multicolumn{3}{c||}{}&\multicolumn{3}{|c||}{\begin{tabular}{c}Iterations
  per\\ run\\
  \end{tabular}}&\multicolumn{2}{c||}{\begin{tabular}{c}CPU per run,\\
  sec\\
 \end{tabular}}&\begin{tabular}{c}CPU per\\
 iteration, sec\\
 \end{tabular}\\
 \hline
\# &Image size&Runs&min&mean&max&mean&max&mean\\
 \hline\hline
A&256$\times$256&6&4&9.00&12&83.4&148.7&8.3\\
\hline
B&512$\times$512&9&4&7.89&11&212.9&318.2&25.9\\
\hline
C&256$\times$256&6&17&17.17&18&189.7&214.7&10.3\\
\hline
D&512$\times$512& 6&16&16.00&16&615.9&768.3&36.0\\
 \hline\hline
 \end{tabular}}
 \caption{\label{tableTV2} Performance of CO{\CG}M; platform: T410 Lenovo laptop, {Intel Core i7 M620 CPU@2.67GHz, 8GB RAM}. Flow solver: interior point method {\tt mosekopt  6.0} \cite{Mosek}}
 \end{center}
 \end{table}

 \bibliographystyle{abbrv}

 \section{Appendix}

\subsection{Proof of Theorem \ref{propredgrad}}
\label{sec:proof p1}
Define
$$
\epsilon_t=f(x_t)-f_*,\quad  \Delta_t=\max\limits_{x\in X}\langle f'(x_t),x_t-x\rangle=\langle f'(x_t),x_t-x_t^+\rangle
$$
where $x^+_t=x_X[f'(x_t)]$. Denoting by $x_*$ an optimal solution to (\ref{problemf}) and invoking the definition of $x_t^+$ and convexity of $f$,
we have
\begin{equation}\label{begin}
\langle f'(x_t),x^+_t-x_t\rangle\le\langle f'(x_t),x_*-x_t\rangle
\le f_*-f(x_t).
\end{equation}
Observing that for a generic GC algorithm we have $f(x_{t+1})\leq f(x_t+\gamma_t(x_t^+-x_t))$  and invoking (\ref{suchthat}), we have
\begin{equation}\label{begin1}
f(x_{t+1})\le f(x_t)+\gamma_t\langle f'(x_t),x_t^+-x_t\rangle +{L\over 2} \gamma_t^2\|x^+_t-x_t\|_X^2
\le f(x_t)-\gamma_t(f(x_t)-f_*) + \half L\gamma_t^2,
\end{equation}
where the concluding $\leq$ is due to (\ref{begin}).
It follows that $\epsilon_{t+1}\leq(1-\gamma_t)\epsilon_t+\half L\gamma_t^2$, whence
\bse
\epsilon_{t+1}&\le&\epsilon_1\prod_{i=1}^t (1-\gamma_i)
+
\half L\sum_{i=1}^t  \gamma_i^2\prod_{k=i+1}^t (1-\gamma_k)\\
&=&
2L\sum_{i=1}^t  (i+1)^{-2}\prod_{k=i+1}^t (1-{2\over k+1}),
\ese
{where, by convention, $\prod_{k=t+1}^t=1$.} Noting that
$
\prod_{k=i+1}^t(1-{2\over k+1})=\prod_{k=i+1}^t {k-1\over k+1}=
{i(i+1)\over t(t+1)},\;\;\;i=1,...,t,
$
 we get
\be
\epsilon_{t+1}\le
2L\sum_{i=1}^t  {i(i+1)\over (i+1)^{2}t(t+1)}\le 
{2L t\over (t+1)^2}\le
{2L (t+2)^{-1}},
\ee{aille+}
what is \rf{then}.
\par To prove \rf{certif1}, observe that setting $\bar{\Delta}_t=\min_{1\leq k\leq t}\Delta_k$,
and invoking (\ref{eq:lowbnd}), (\ref{eq:lowbnd1}) we clearly have
$$
\begin{array}{rcl}
f(\bar{x}_t)-f_t^*&=&\min_{1\leq k\leq t}[f(\bar{x}_t)-f_{*,k}]=\min_{1\leq k\leq t}[f(\bar{x}_t)-f(x_k)+\Delta_k]\\
&\leq& \min_{1\leq k\leq t}\Delta_k =\bar{\Delta}_t\\
\end{array}
$$
(we have used the fact that $f(\bar{x}_t)\leq f(x_k)$, $k\leq t$, by definition of $\bar{x}_t$). We see that in order to prove (\ref{certif1}), it suffices to prove that
\begin{equation}\label{certif11}
\bar{\Delta}_t\leq {4.5 L\over t-2},\,\,t=5,6,...
\end{equation}
To verify (\ref{certif11}), note that by the first inequality in (\ref{begin1})
\begin{equation}\label{moreofthesame}
\gamma_k\Delta_k\equiv \gamma_k\langle f'(x_k),x_k-x^+_k\rangle\le \epsilon_k-\epsilon_{k+1}+\half L\gamma_k^2.
\end{equation}
Assuming $t>2$ and summing up these inequalities over $k$ varying from $t_0=\rfloor t/2\lfloor$ to $t$ (here $\rfloor a\lfloor$ stands for the largest integer strictly smaller than $a$), we obtain
\[
\sum_{k=t_0}^t\gamma_k\Delta_k\le \epsilon_{t_0}+\half L\sum_{k=t_0}^t\gamma_k^2,
\]
 and therefore
\be
\bar{\Delta}_t\sum_{k=t_0}^t\gamma_k\le \epsilon_{t_0}+\half L\sum_{k=t_0}^t\gamma_k^2.
\ee{aille}
Observe that $
\sum_{k=t_0}^t\gamma_k=2\sum_{k=t_0}^t(k+1)^{-1}\ge 2[\ln(t+1)-\ln(t_0+1)]\ge 2\ln(2)
$
 and
 $
 \sum_{k=t_0}^t\gamma_k^2=4\sum_{k=t_0}^t (k+1)^{-2}\le 4[t_0^{-1}-t^{-1}]\le {4(t+2)\over t(t-2)}.
$
Assuming $t>4$ (so that $t_0\geq2$) and substituting into \rf{aille} the bound \rf{then} for $\epsilon_{t_0}$ we obtain
\[
\bar{\Delta}_t\le {2L\over t}+{L(t+2)\over t(t-2)\ln(2)}\le 4.5 L(t-2)^{-1},
\]
as required in (\ref{certif11}).
\qed

\subsection{Proof of Theorem \ref{TheMain2}}
The proof, up to minor modifications, goes back to \cite{Lem:Nem:Nes:1995}, see also \cite{Nesterov:2004,NemJud1}; we provide it here to make the paper self-contained. W.l.o.g. we can assume that we are in the nontrivial case (see description of the algorithm).
\par
{\bf 1$^0$.} As it was explained when describing the method, whenever stage $s$ takes place, we have $[0<]\rho_1\leq\rho_s\leq\rho_*$, and $\rho_{s-1}<\rho_s$, provided $s>1$. Therefore by the termination rule, the output $\bar{\rho}$, $\bar{x}$ of the algorithm, if any, satisfies $\bar{\rho}\leq\rho_*$, $f(\bar{x})\leq\epsilon$. Thus, (i) holds true, provided that the algorithm does terminate. Thus, all we need is to verify (ii) and (iii). \par
{\bf 2$^0$.} Let us prove (ii). Let $s\geq1$ be such that stage $s$ takes place. Setting $X=K[\rho_s]$, observe that $X-X\subset \{x\in E:\|x\|\leq 2\rho_s\}$, whence $\|\cdot\|\leq 2\rho_s\|\cdot\|_X$, and therefore the relation (\ref{suchthatinitial}) implies the validity of (\ref{suchthat}) with $L=4\rho_s^2L_f$. Now, if stage $s$ does not terminate in course of some number $t$ steps, then, in the notation from the description of the algorithm, $f(\bar{x}_t)>\epsilon$ and $f_*^t<{3\over 4}f(\bar{x}_t)$, whence $f(\bar{x}_t)-f_*^t> \epsilon/4$. By Theorem \ref{propredgrad}.ii, the latter is possible only when $4.5L/(t-2)>\epsilon/4$. Thus,
$t\leq \max\left[5,2+{72\rho_s^2L_f\over\epsilon}\right]$. Taking into account that $\rho_s\leq\rho_*$, (ii) follows.
\par
{\bf 3$^0$.} Let us prove (iii). This statement is trivially true when the number of stages is 1. Assuming that it is not the case, let $S\geq1$ be such that the stage $S+1$ takes place. For every $s=1,...,S$, let $t_s$ be the last step of stage $s$, and let $u_s$, $\ell^s(\cdot)$ be what in the notation from the description of stage $s$ was denoted $f(\bar{x}_{t_s})$ and $\ell^{t_s}(\rho)$. Thus,
$u_s>\epsilon$ is an upper bound on $\Opt(\rho_s)$, $\ell_s:=\ell^s(\rho_s)$ is a lower bound on $\Opt(\rho_s)$ satisfying $\ell_s\geq 3u_s/4$, and $\ell^s(\cdot)$ is a piecewise linear convex in $\rho$ lower bound on $\Opt(\rho)$, $\rho\geq0$, and $\rho_{s+1}>\rho_s$ is the smallest positive root  of $\ell^s(\cdot)$.  Let also $-g_s$ be a subgradient of $\ell^s(\cdot)$ at $\rho_s$. Note that $g_s>0$ due to $\rho_{s+1}>\rho_s$ combined with $\ell^s(\rho_s)>0$, $\ell^s(\rho_{s+1})=0$, and by the same reasons combined with convexity of $\ell^s(\cdot)$ we have
\begin{equation}\label{rhonext}
\rho_{s+1}-\rho_s\geq \ell_s/g_s,
\end{equation}
and, as we have seen,
\begin{equation}\label{then12}
1\leq s\leq S\Rightarrow \left\{
\begin{array}
{ll}
(a)&u_s>\epsilon,\\
(b)&u_s\geq\Opt(\rho_s)\geq \ell_s\geq {3\over 4}u_s,\\
(c)&\ell_s-g_s(\rho-\rho_s) \leq\Opt(\rho),\,\rho\geq 0.\\
\end{array}\right..
\end{equation}
Assuming $1<s\leq S$ and applying  (\ref{rhonext}), we get $\rho_s-\rho_{s-1}\geq {3\over 4}u_{s-1}/g_{s-1}$, whence, invoking (\ref{then12}),
$$u_{s-1}\geq \Opt(\rho_{s-1})\geq \ell_s+g_s[\rho_{s-1}-\rho_s]\geq {3\over 4}u_s+{3\over 4}u_{s-1}{g_s\over g_{s-1}}.$$
The resulting inequality implies that
${u_s\over u_{s-1}}+{g_s\over g_{s-1}}\leq {4\over 3}$,  whence ${u_sg_s\over u_{s-1}g_{s-1}}\leq (1/4)(4/3)^2=4/9$.
It follows that
\begin{equation}\label{thatso}
\sqrt{u_sg_s}\leq (2/3)^{s-1}\sqrt{u_1g_1}, \,\,1\leq s\leq S.
\end{equation}
Now, since the first iterate of the first stage is $0$, we have $u_1\leq f(0)$, while (\ref{then12}) applied with $s=1$ implies that $f(0)=\Opt(0)\geq \ell_1+\rho_1g_1\geq \rho_1g_1$, whence $u_1g_1\leq f(0)/\rho_1=d$. Further, by (\ref{rhonext}) $g_s\geq\ell_s/(\rho_{s+1}-\rho_s)\ge
q \ell_s/\rho_*\geq {3\over 4} u_s/\rho_*$, where the concluding inequality is given by (\ref{then12}). We see that $u_sg_s\geq {3\over 4}u_s^2/\rho_*\geq{3\over 4}\epsilon^2/\rho_*$. This lower bound on $u_sg_s$ combines with the bound $u_1g_1\leq d$ and with (\ref{thatso}) to imply that
$$
\epsilon\leq \sqrt{4/3}(2/3)^{s-1}\sqrt{d\rho_*},\,1\leq s\leq S.
$$
Finally observe that by the definition of $\rho_*$ and due to the fact that $\|x[f'(0)]\|=1$ in the nontrivial case, we have
$$
0\leq f(\rho_*x[f'(0)])\leq f(0)+\rho_*\langle f'(0),x[f'(0)]\rangle +{1\over 2}L_f\rho_*^2=f(0)-\rho_*d+\half L_f\rho_*^2
$$
(we have used (\ref{suchthatinitial}) and the definition of $d$), whence $\rho_*d\leq f(0)+\half L_f\rho_*^2$ and therefore
$$
\epsilon\leq \sqrt{3/4}(2/3)^{s-1}\sqrt{f(0)+\half L_f\rho_*^2},\,1\leq s\leq S.
$$
Since this relation holds true for every $S\geq1$ such that the stage $S+1$ takes place, (iii) follows. \qed

\subsection{Proof of Theorem \ref{pro:cgco}}
By definition of $z_t$ we have $z_t\in K^+$ for all $t$ and $F(0)=F(z_1)\ge F(z_2)\ge ...$,
whence $r_t\leq D_*$ for all $t$ by Assumption A. Besides this, $r_*\leq D_*$ as well. Let now $\epsilon_t=F(z_t)-F_*$,  $z_t=[x_t;r_t]$, and let
$z^+_t=[x^+_t,r^+_t]$  be a minimizer, as given by Lemma \ref{lem:evident},  of the linear form $\langle F'(z_t),z\rangle$ of $z\in E^+$ over the set $K^+[r_*]=\{[x;r]: x\in K,\|x\|\leq r\leq r_*\}$. Recalling that $F'(z_t)=[f'(x_t);\kappa]$ and that $r_t\leq D_*\leq D^+$, Lemma \ref{lem:evident} implies that $z^+_t\in \Delta(z_t)$.  By definition of $z_t^+$ and convexity of $F$ we have
$$
\begin{array}{rcl}
\langle [f'(x_t);\kappa],z_t-z_t^+\rangle&=&\langle f'(x_t),x_t-x_t^+\rangle +\kappa (r_t-r_t^+)\\
&\geq& \langle f'(x_t),x_t-x_*\rangle +\kappa (r_t-r_*)\\
&=&
 \langle F'(z_t),\,z_t-z_*\rangle
\ge F(z_t)-F(z_*)
 =\epsilon_t.\\
 \end{array}
$$
Invoking (\ref{suchthat}), it follows that for $0\leq s\leq 1$ one has
\bse
F(z_t+s(z_t^+-z_t))&\leq& F(z_i)+s\langle [f'(x_t);\kappa],z_t^+-z_t\rangle +{L_f s^2\over 2}\|x(z_t^+)-x(z_t)\|^2\\
&\leq&
F(z_t)-s\epsilon_t
+\half L_f s^2(r_t+D_*)^2
\ese
{using that $\|x(z_t^+)\|\leq r_t^+$ and $\|x(z^t)\|\leq r_t$
due to $z_t^+,z_t\in K^+$, and that $r_t^+\leq r_*\leq D_*$.}
By (\ref{cgco}) we have
$$
F(z_{t+1})\leq\min_{0\leq s\leq 1}F(z_t+s(z_t^+-z_t))\leq F(z_t)+\min\limits_{0\leq s\leq1} \left\{-s\epsilon_t
+\half L_f s^2(r_t+D_*)^2\right\},
$$
and we arrive at the recurrence
\be
\epsilon_{t+1}\leq\epsilon_t-\left\{\begin{array}{ll}{\epsilon_t^2\over 2L_f (r_t+D_*)^2},&\epsilon_t\leq L_f (r_t+D_*)^2\\
 \epsilon_t-\half L_f (r_t+D_*)^2,&\epsilon_t> L_f (r_t+D_*)^2\\
\end{array}\right., t=1,2,...
\ee{recursive}
When $t=1$, this recurrence, in view of $z_1=0$, implies that $\epsilon_2\leq \half L_fD_*^2$.
Let us show by induction in $t\geq2$ that
\begin{equation}\label{induction}
\epsilon_{t}\leq \bar{\epsilon}_t:={8L_f D_*^2\over t+14},\,t=2,3,...
\end{equation}
thus completing the proof. We have already seen that (\ref{induction}) is valid for $t=2$. Assuming that (\ref{induction}) holds true for $t=k\geq2$, we have $\epsilon_{k}\leq \half L_f D_*^2$ and therefore
$\epsilon_{k+1}\leq \epsilon_{k}-{1\over 8L_f D_*^2}\epsilon_{k}^2$
by \rf{recursive} combined with $0\leq r_k\leq D_*$. Now, the function $ s-{1\over 8L_f D_*^2}s^2$ is nondecreasing on the segment $0\leq s\leq 4L_f D_*^2$ which contains $\bar{\epsilon}_k$ and $\epsilon_k\leq\bar{\epsilon}_k$, whence
\bse
\epsilon_{k+1}&\leq& \epsilon_{k}-{1\over 8L_f D_*^2}\epsilon_{k}^2\leq \bar{\epsilon}_{k}-{1\over 8L_f D_*^2}\bar{\epsilon}_{k}^2
= \left[{8L_f D_*^2\over k+14}\right]-{1\over 8L_f D_*^2}\left[{8L_f D_*^2\over k+14}\right]^2\\
&=& {8L_f D_*^2(k+13)\over (k+14)^2}\leq {8L_f D_*^2\over (k+1)+14},\\
\ese
so that (\ref{induction}) holds true for $t=k+1$.  \qed
\subsection{Proofs for Section \ref{sec:examples}}
\paragraph{Proof of Lemma \ref{lemflow}.}  As we have already explained, (\ref{network}) is solvable, so that $z$ is well defined. Denoting by $(s^*,r^*)$ an optimal solution to (\ref{network}) produced, along with $z$, by our solver, note that the characteristic property of $z$ is the relation
$$
(s^*,r^*)\in\Argmax_{s,r}\{s+\langle z, Pr-s\eta\rangle: 0\leq r\leq \mathbf{e}\}.
$$
Since the column sums in $P$ are zeros and $\eta$ is with zero sum of entries, the above characteristic property of $z$ is preserved when passing from $z$ to $\bar{z}$, so that we may assume from the very beginning that $z=\bar{z}$ is a zero mean image. Now, $P=[Q,-Q]$, where $Q$ is the incidence matrix of the network obtained from $G$ by eliminating backward arcs. Representing a flow $r$ as $[r_f;r_b]$, where the blocks are comprised, respectively, of flows in the forward and backward arcs, and passing from $r$ to $\rho=r_f-r_b$, our characteristic property of $z$ clearly implies the relation
$$
(s^*,\;\rho^*:=r^*_f-r^*_b)\in\Argmax_{s,\rho}\{s+\langle z,Q\rho-s\eta\rangle:\|\rho\|_\infty\leq1\}.
$$
It follows that
\begin{equation}\label{lagrange}
\begin{array}{ll}
(a)&\langle z,\eta\rangle=1,\\
(b)&\|\rho^*\|_\infty\leq1,\\
(c)&(Q^*z)_\gamma=\left\{\begin{array}{ll}\leq0, &\rho^*_\gamma=-1,\\
=0,&\rho^*_\gamma\in(-1,1),\\
\geq0,&\rho^*_\gamma=1,\\
\end{array}\right.\hbox{\ for all forward arcs $\gamma$},\\
(d)&Q\rho^*=s^*\eta.\\
\end{array}
\end{equation}
(\ref{lagrange}.$d$) and (\ref{lagrange}.$a$) imply that $\langle Q^*z,\rho^*\rangle =s^*$, while (\ref{lagrange}.$c$) says that $\langle Q^*z,\rho^*\rangle = \|Q^*z\|_1$, and
$s^*=\|Q^*z\|_1$.  By (\ref{lagrange}.$a$) $z\neq0$, and thus $z$ is a nonzero image with zero mean; recalling what $Q$ is, the first $n(n-1)$ entries in $Q^*z$ form $\nabla_i z$, and the last $n(n-1)$ entries form $\nabla_jz$, so that $\|Q^*z\|_1=\TV(z)$. The gradient field of a nonzero image with zero mean cannot be identically zero, whence $\TV(z)=\|Q^*z\|_1=s^*>0$. Thus $x[\eta]=-z/\TV(z)=-z/s^*$ is well defined and $\TV(x[\eta])=1$, while by (\ref{lagrange}.$a$) we have $\langle x[\eta],\eta\rangle=-1/s^*$. Finally, let $x\in\CTV$, implying that $Q^*x$ is the concatenation of $\nabla_ix$ and $\nabla_jx$ and thus $\|Q^*x\|_1=\TV(x)\leq1$. Invoking (\ref{lagrange}.$b,d$), we get $-1\leq \langle Q^*x,\rho^*\rangle=\langle x, Q\rho^*\rangle =s^*\langle x,\eta\rangle$, whence $\langle x,\eta\rangle \geq -1/s^*=\langle x[\eta],\eta\rangle$, meaning that $x[\eta]\in\CTV$ is a minimizer of $\langle \eta,x\rangle$ over $x\in\CTV$. \qed

\paragraph{Proof of Proposition \ref{laplace}.} In the sequel, for a real-valued function $x$ defined on a finite set (e.g., for an image), $\|x\|_p$ stands for the $L_p$ norm of the function corresponding to the counting measure on the set (the mass of every point from the set is 1). Let us fix $n$ and $x\in M^n_0$ with $\TV(x)\leq1$; we want to prove that
\begin{equation}\label{target}
\|x\|_2\leq \cC\sqrt{\ln(n)}
\end{equation}
with appropriately selected {\sl  absolute constant} $\cC$.
\par
{\bf 1$^0$.} Let $\oplus$ stand for addition, and $\ominus$ -- for substraction of integers modulo $n$; $p\oplus q=(p+q) \,\hbox{mod} \,n\in\{0,1,...,n-1\}$ and similarly for $p\ominus q$. Along with discrete partial derivatives $\nabla_i x$, $\nabla_jx$, let us define their periodic versions $\widehat{\nabla}_ix$, $\widehat{\nabla}_jx$:
$$
\widehat{\nabla}_ix(i,j)= x(i\oplus 1,j)-x(i,j):\Gamma_{n,n}\to \bR,\,\,
\widehat{\nabla}_jx(i,j)= x(i,j\oplus1)-x(i,j):\Gamma_{n,n}\to \bR,
$$
same as periodic Laplacian $\widehat{\Delta} x$:
$$
\widehat{\Delta} x=x(i,j)-{1\over 4}\left[x(i\ominus1,j)+x(i\oplus1,j)+x(i,j\ominus1)+x(i,j\oplus1)\right]:\Gamma_{n,n}\to\bR.
$$
For every $j$, $0\leq j<n$, we have $\sum_{i=0}^{n-1} \widehat{\nabla}_ix(i,j)=0$ and $\nabla_ix(i,j)=\widehat{\nabla}_ix(i,j)$ for $0\leq i<n-1$, whence $\sum_{i=0}^{n-1}|\widehat{\nabla}_i(x)|\leq 2\sum_{i=0}^{n-1}|\nabla_ix(i,j)|$ for every $j$, and thus $\|\widehat{\nabla}_ix\|_1\leq 2\|\nabla_ix\|_1$. Similarly, $\|\widehat{\nabla}_jx\|_1\leq2\|\nabla_jx\|_1$, and we conclude that
\begin{equation}\label{leq2}
\|\widehat{\nabla}_ix\|_1+\|\widehat{\nabla}_jx\|_1\leq 2.
\end{equation}
\par
{\bf 2$^0$.} Now observe that for $0\leq i,j<n$ we have
$$
\begin{array}{rcl}
x(i,j)&=&x(i\ominus1,j)+\widehat{\nabla}_ix(i\ominus1,j)\\
x(i,j)&=&x(i\oplus1,j)-\widehat{\nabla}_ix(i,j)\\
x(i,j)&=&x(i,j\ominus1)+\widehat{\nabla}_jx(i,j\ominus1)\\
x(i,j)&=&x(i,j\oplus1)-\widehat{\nabla}_jx(i,j)\\
\end{array}
$$
whence
\begin{equation}\label{Delta}
\widehat{\Delta} x(i,j)={1\over4}\left[\widehat{\nabla}_ix(i\ominus1,j)-\widehat{\nabla}_ix(i,j)+\widehat{\nabla}_jx(i,j\ominus1)-
\widehat{\nabla}_jx(i,j)\right]
\end{equation}
Now consider the following linear mapping from $M^n\times M^n$ into $M^n$:
\begin{equation}\label{B}
B[g,h](i,j)={1\over4}\left[g(i\ominus1,j)-g(i,j)+h(i,j\ominus1)-h(i,j)\right],\,[i;j]\in\Gamma_{n,n}.
\end{equation}
From this definition and (\ref{Delta}) it follows that
\begin{equation}\label{Deltaxiseq}
\widehat{\Delta}x=B[\widehat{\nabla}_ix,\widehat{\nabla}_jx].
\end{equation}
\par
{\bf 3$^0$.} Observe that $B[g,h]$ always is an image with zero mean. Further, passing from images $u\in M^n$ to their 2D Discrete Fourier Transforms $\DFT[u]$:
$$
\DFT[u](p,q)=\sum_{0\leq r,s<n} u(r,s)\exp\{-2\pi\imath(pr+qs)/n\},\,[p;q]\in\Gamma_{n,n},
$$
we immediately see that every image $u$ with zero mean is the periodic Laplacian of another, uniquely, defined, image $X[u]$ with zero mean, with $X[u]$ given by
\begin{equation}\label{capitals}
\begin{array}{c}
\DFT[X[u]](p,q)=Y[u](p,q):=\left\{\begin{array}{ll}0,&p=q=0\\
{\dft[u](p,q)\over D(p,q)}, &0\neq [p;q]\in \Gamma_{n,n}\\
\end{array}\right.,\,[p;q]\in \Gamma_{n,n},\\
D(p,q)=1-{1\over 2}[\cos(2\pi p/n)+\cos(2\pi q/n)],\,\,[p;q]\in\Gamma_{n,n}.
\end{array}
\end{equation}
In particular, invoking (\ref{Deltaxiseq}), we get
$$
\DFT[x]=Y[B[\widehat{\nabla}_ix,\widehat{\nabla}_jx]].
$$
By Parseval identity, $\|\DFT[x]\|_2=n\|x\|_2$, whence
$$
\|x\|_2=n^{-1}\|Y[B[\widehat{\nabla}_ix,\widehat{\nabla}_jx]]\|_2.
$$
Combining this observation with (\ref{leq2}), we see that in order to prove (\ref{target}), it suffices to check that
\begin{quote}
(!) {\sl Whenever $g,h\in M^n$ are such that \[(g,h)\in G:=\{(g,h)\in M^n\times M^n: \|g\|_1+\|h\|_1\leq2\},\] we have}
\begin{equation}\label{exclamation}
\|Y[B[g,h]]\|_2\leq n\cC\sqrt{\ln(n)}.
\end{equation}
\end{quote}
\par
{\bf 4$^0$.} A good news about (!) is that since $Y[B[g,h]]$ is linear in $(g,h)$, in order to justify (!), it suffices to prove that (\ref{exclamation}) holds true for the extreme point of $G$, i.e., (a) for pairs where $h\equiv 0$ and $g$ is an image which is equal to 2 at some point of $\Gamma_{n,n}$ and vanishes outside of this point, and (b) for pairs where $g\equiv 0$ and $h$ is an image which is equal to 2 at some point of $\Gamma_{n,n}$ and vanishes outside of this point. Task (b) clearly reduces to task (a) by swapping the coordinates $i,j$ of points from $\Gamma_{n,n}$, so that we may focus solely on task (a). Thus, assume that $g$ is a cyclic shift of the image $2\delta$:
$$
g(i,j)\equiv 2\delta(i\ominus r,j\ominus s),\,\,\delta(i,j)=\left\{\begin{array}{ll} 1,&[i;j]=[0;0]\\
0,&[i;j]\neq[0;0]\\
\end{array}\right.,\,[i;j]\in\Gamma_{n,n}.
$$
From (\ref{B}) it follows that then $B[g,0]$ is a cyclic shift of $B[2\delta,0]$, whence $|\DFT[B[g,0]](p,q)|=|\DFT[B[2\delta,0]](p,q)|$ for all $[p;q]\in\Gamma_{n,n}$, which, by (\ref{capitals}), implies that $|Y[B[g,0]](p,q)|=|Y[B[2\delta,0]](p,q)|$ for all $[p;q]\in\Gamma_{n,n}$. The bottom line is that {\sl all we need is to verify that {\rm (\ref{exclamation})} holds true for $g=2\delta,h=0$,} or, which is the same, that with
\begin{equation}\label{eqy}
y(p,q)={(1-\exp\{2\pi\imath p/n\})\over 2[1-{1\over 2}[\cos(2\pi p/n)+\cos(2\pi q/n)]]}
\end{equation}
where the right hand side by definition is $0$ at $p=q=0$, it holds
$$
C_n:=\sum_{p,q=0}^{n-1} |y(p,q)|^2 \leq n^2\cC^2\ln(n).
$$
Now, (\ref{eqy}) makes sense for all $[p;q]\in\bZ^2$ (provided that we define the right hand side as zero at all points of $\bZ^2$ where the denominator in (\ref{eqy}) vanishes, that is, at all point where $p,q$ are integer multiples of $n$) and defines $y$ as a double-periodic, with periods $n$ in $p$ and in $q$, function of $[p;q]$. Therefore, setting $m=\hbox{Floor}(n/2)\geq 1$ and $W=\{[p;q]\in\bZ^2: -m\leq p,q<n-m\}$, we have
$$
C_n=\sum_{0\neq [p;q]\in W} |y(p,q)|^2=\sum_{[p;q]\in W} {|1-\exp\{2\pi\imath p/n\}|^2\over 4|1-{1\over 2}[\cos(2\pi p/n)+\cos(2\pi q/n)]|^2}.
$$
Setting $\rho(p,q)=\sqrt{p^2+q^2}$, observe that when $0\neq[p;q]\in W$, we have $|1-\exp\{2\pi\imath p/n\}|\leq C_1n^{-1}\rho(p,q)$ and
$2[1-{1\over 2}[\cos(2\pi p/n)+\cos(2\pi q/n)]]\geq C_2n^{-2}\rho^2(p,q)$ with positive absolute constants $C_1,C_2$, whence
$$
C_n\leq(C_1/C_2)^2\sum_{0\neq [p;q]\in W} n^2\rho^{-2}(p,q).
$$
With appropriately selected absolute constant $C_3$ we have
$$
\sum_{0\neq [p;q]\in W} \rho^{-2}(p,q)\leq C_3\int_{1}^n r^{-2}rdr=C_3\ln(n).
$$
Thus, $C_n\leq (C_1/C_2)^2C_3n^2\ln(n)$, meaning that (\ref{exclamation}), and thus (\ref{target}), holds true with $\cC=\sqrt{C_3}C_1/C_2$. \qed
\end{document}